%% file: manuscript-amspreprint.tex
\pdfoutput=1
\IfFileExists{./prepreamble-amspreprint.sty}{\RequirePackage[packages,theorems]{./prepreamble-amspreprint}}{}

\documentclass[]{amsart}

\RequirePackage{biblatex}
\RequirePackage{ltxcmds}
\IfFileExists{./preamble-amspreprint.sty}{\RequirePackage[packages,theorems]{./preamble-amspreprint}}{}

\usepackage{manuscript}
\input{macros}


\usepackage[ruled, linesnumbered,algo2e]{algorithm2e}
\usepackage{amsfonts}
\usepackage{mathtools}
\usepackage{amsthm}
\usepackage[normalem]{ulem}
\usepackage{subcaption}
\usepackage{amsmath,amssymb,bm}
\usepackage{enumerate}
\usepackage{hyperref} 
\usepackage{cleveref}
\usepackage{algorithm, algpseudocode}
\usepackage{psfrag,epsf,color}
\usepackage{booktabs}
\usepackage[dvipdfm]{graphicx} 
\usepackage{bmpsize}
\usepackage[utf8]{inputenc}
\DeclareUnicodeCharacter{2082}{\textsubscript{2}}

\IfFileExists{./postpreamble-amspreprint.sty}{\RequirePackage[packages,theorems]{./postpreamble-amspreprint}}{}

\makeatletter
\@ifpackageloaded{hyperref}{%
	\hypersetup{
		pdftitle = {Fast solution of Sylvester-structured systems for spatial source separation of the Cosmic Microwave Background\protect},
		pdfauthor = {Dung Pham, Kirk M. Soodhalter, Simon Wilson},
		pdfkeywords = {Cosmic Microwave Background, Gaussian field, Krylov subspace, source separation, Conjugate Gradients, Sylvester Equations},
		pdfcreator = {Created using the Scoop Template Engine version 1.2.0.}
	}
}{
	\pdfinfo{
		/Title (Fast solution of Sylvester-structured systems for spatial source separation of the Cosmic Microwave Background\protect)
		/Author (Dung Pham, Kirk M. Soodhalter, Simon Wilson)
		/Subject ()
		/Keywords (Cosmic Microwave Background, Gaussian field, Krylov subspace, source separation, Conjugate Gradients, Sylvester Equations)
		/Creator (Created using the Scoop Template Engine version 1.2.0.)
	}
}
\makeatother

\title[FAST SOLVERS FOR SYLVESTER-STRUCTURED CMB PROBLEMS]{Fast solution of Sylvester-structured systems for spatial source separation of the Cosmic Microwave Background\protect}

\author{Dung Pham}
\address[D. Pham]{School of Computer Science and Statistics Trinity College Dublin College Green, Dublin 2 Ireland}
\email{phamd@tcd.ie}

\author{Kirk M. Soodhalter}
\address[K. M. Soodhalter]{School of Mathematics Trinity College Dublin College Green, Dublin 2 Ireland}
\email{ksoodha@maths.tcd.ie}
\urladdr{https://math.soodhalter.com}

\author{Simon Wilson}
\address[S. Wilson]{School of Computer Science and Statistics Trinity College Dublin College Green, Dublin 2 Ireland}
\email{SWILSON@tcd.ie}

\date{\today}

\dedicatory{}

\begin{document}

% Insert the abstract.
\begin{abstract}
\input{abstract.tex}
\end{abstract}

% Insert the keywords.

% Insert the Mathematics Subject Classification.
\makeatletter
\ltx@ifpackageloaded{hyperref}{%
\subjclass[2010]{\href{https://mathscinet.ams.org/msc/msc2020.html?t=65F10}{65F10}, \href{https://mathscinet.ams.org/msc/msc2020.html?t=65F50}{65F50}, \href{https://mathscinet.ams.org/msc/msc2020.html?t=65F08}{65F08}}
}{%
}
\makeatother

% Typeset the opening page.
\maketitle

% Insert the document body.
\input{content.tex}

% Insert the appendix.
\appendix
\input{appendix.tex}

% Insert the bibliography.
\printbibliography

\end{document}

%% file: macros.tex
\def\bcalY{\boldsymbol{\cal Y}}
\def\bcalS{\boldsymbol{\cal S}}
\def\bcalE{\boldsymbol{\cal E}}
\def\beps{\boldsymbol\epsilon}
\def\bmu{\boldsymbol\mu}
\def\cal{\mathcal}

\def\nm#1{\left\|#1\right\|}
\def\norm#1{\nm{#1}_{2}}

\def\R{\mathbb{R}}

\def\K{\mathbb{K}}

\def\Rn{\R^n}

\def\Rmm{\R^{m\times m}}
\def\Rnn{\R^{n\times n}}

\def\BA{{\bf A}}

  \def\CK{{\cal K}}

  \def\CO{{\cal O}}
  \def\CP{{\cal P}}

  \def\CV{{\cal V}}

\def\BAs{\BA{\kern-1.5pt}}

\def\CPs{\CP{\kern-0.8pt}}

\mathcode`\@="8000
{\catcode`\@=\active \gdef@{\mkern1mu}}
\def\mydate{\number\day\ {\ifcase\month \or January\or February\or
              March\or April\or May\or June\or July\or August\or
              September\or October\or November\or December\fi}
\number\year}
\def\vek#1{\mathbf{#1}}

\providecommand{\prn}[1]{\left(#1\right)}

\def\curl#1{\left\{#1\right\}}

\newcommand\restr[2]{{% we make the whole thing an ordinary symbol
  \left.\kern-\nulldelimiterspace % automatically resize the bar with \right
  #1 % the function
  \vphantom{\big|} % pretend it's a little taller at normal size
  \right|_{#2} % this is the delimiter
  }}

\def\be{\begin{equation}}
\def\ee{\end{equation}}
\def\bea{\begin{eqnarray}}
\def\eea{\end{eqnarray}}
\def\nn{\nonumber}
\def\mand{\mbox{\ \ \ and\ \ \ }}

\def\mwith{\mbox{\ \ \ with\ \ \ }}
\def\mwhere{\mbox{\ \ \ where\ \ \ }}

\def\mif{\mbox{\ \ \ if\ \ \ }}
\def\mthen{\mbox{\ \ \ then\ \ \ }}

\def\bbmat{\begin{bmatrix}}
\def\ebmat{\end{bmatrix}}

\def\balg{\begin{algorithm}}
\def\ealg{\end{algorithm}}

\def\bthm{\begin{theorem}}
\def\ethm{\end{theorem}}
\def\blem{\begin{lemma}}
\def\elem{\end{lemma}}
\def\bprop{\begin{proposition}}
\def\eprop{\end{proposition}}
\def\bcor{\begin{corollary}}
\def\ecor{\end{corollary}}
\def\bdefin{\begin{definition}}
\def\edefin{\end{definition}}
\def\bc{\begin{cases}}
\def\ec{\end{cases}}
\newcommand\bproof[1]{\par\addvspace{1ex} \indent\textit{Proof.}\ \ #1}
\def\eproof{\hfill\cvd\linebreak\indent}
\newtheorem{exple}{Example}
\def\bex{\begin{exple}}
\def\eex{\end{exple}}
\def\bexer{\begin{exercise}}
\def\eexer{\end{exercise}}
\def\bconj{\begin{conjecture}}
\def\econj{\end{conjecture}}
\def\bass{\begin{assumption}}
\def\eass{\end{assumption}}
\def\bnot{\begin{notation}}
\def\enot{\end{notation}}
\def\brem{\begin{remark}}
\def\erem{\end{remark}}
\def\diag{{\rm diag}}

\def\blankfootnote#1{\let\thefootnote\relax\footnotetext{#1}}
\def\cvd{~\vbox{\hrule\hbox{%
  \vrule height1.3ex\hskip0.8ex\vrule}\hrule } }

\def\bitem{\begin{item}}
\def\eitem{\end{itemize}}
\def\benum{\begin{enumerate}}
\def\eenum{\end{enumerate}}
  
\DeclareMathSymbol{\dprod}{\mathbin}{operators}{"3A}

\providecommand{\matlab}{\textsc{Matlab}}
\providecommand{\file}[1]{\texttt{\nolinkurl{#1}}}

%% file: abstract.tex
Implementation of many statistical methods for large, multivariate data sets requires one to solve a linear system that, depending on the method, is of the dimension of the number of observations or each individual data vector.  This is often the limiting factor in scaling the method with data size and complexity. In this paper we illustrate the use of Krylov subspace methods to address this issue in a statistical solution to a source separation problem in cosmology where the data size is prohibitively large for direct solution of the required system.  Two distinct approaches, adapted from techniques in the literature, are described: one that uses the method of conjugate gradients directly to the Kronecker-structured problem and another that reformulates the system as a Sylvester matrix equation.  We show that both approaches produce an accurate solution within an acceptable computation time and with practical memory requirements for the data size that is currently available.

%% file: content.tex
% %%%%%%%%%%%%%%%%%%%%%%%%%%%%%%%%%%%%%%%%%%%%%%%%%%%%%%%%%%%%%%%%%%%%%%%%%%%%%%%%%%%%%%%%%%%%%%%%%%%%%%%%%%%%%%%%%%%%%%%%%%
% %%%%%%%%%%%%%%%%%%%%%%%%%%%%%%%%%%%%%%%%%%%%%%%%%%%%%%%%%%%%%%%%%%%%%%%%%%%%%%%%%%%%%%%%%%%%%%%%%%%%%%%%%%%%%%%%%%%%%%%%%%
\section{Introduction}
% %%%%%%%%%%%%%%%%%%%%%%%%%%%%%%%%%%%%%%%%%%%%%%%%%%%%%%%%%%%%%%%%%%%%%%%%%%%%%%%%%%%%%%%%%%%%%%%%%%%%%%%%%%%%%%%%%%%%%%%%%%
% %%%%%%%%%%%%%%%%%%%%%%%%%%%%%%%%%%%%%%%%%%%%%%%%%%%%%%%%%%%%%%%%%%%%%%%%%%%%%%%%%%%%%%%%%%%%%%%%%%%%%%%%%%%%%%%%%%%%%%%%%%
In this paper we describe an approach to solving a sparse, high-dimensional linear system that arises in the context of a particular form of the statistical model known as factor analysis \cite{bartholomew11}, also known in the machine learning literature as independent components analysis \cite[Chapter 12]{bishop07}.  This model is used in the analysis of multivariate data, most commonly for:
\begin{itemize}
\item Dimension reduction -- identifying a lower dimensional linear subspace of the data space that still captures a large amount of the data variation, to which further analysis is then applied. This is one approach to handling data whose dimension makes analysis computationally infeasible;
\item Interpretation -- identifying linear combinations of data components that are able to be interpreted in some way, for the purposes of gaining insights into the data generating process or learning about those linear combinations. One application of this is source separation, \emph{the motivating example for this work}.
\end{itemize}
This work is motivated by a problem arising out of data analysis in cosmology.
Our approach makes use of Krylov subspaces.  Methods built on these subspaces have been used to apply statistical methods to scale to larger data sets, starting with \cite{simpson08} and then in a series of papers to \cite{aune14}.  Like this work, those papers focused on the problems associated with working with high-dimensional Gaussian distributions, where statistical inference and prediction usually involve  operations on its variance matrix that is the dimension of the data. This matrix is often designed to be sparse and hence suited to Krylov subspace methods.  In the specific context of CMB reconstruction, Krylov subspace iterative methods employing subspace recycling acceleration have been used to great effectiveness \cite{PapezGrigoriStompor:2020:1}.

This paper is organised as follows.  Section \ref{sec:motivation} discusses the source separation problem for the Cosmic Microwave Background, the motivating example that led to this work.  Section \ref{sec:solution} describes our approaches to solving the linear system that arises out of the example, which are then evaluated in Section \ref{sec:evaluation}.  Section \ref{sec:conclusion} contains some concluding remarks.

% %%%%%%%%%%%%%%%%%%%%%%%%%%%%%%%%%%%%%%%%%%%%%%%%%%%%%%%%%%%%%%%%%%%%%%%%%%%%%%%%%%%%%%%%%%%%%%%%%%%%%%%%%%%%%
% %%%%%%%%%%%%%%%%%%%%%%%%%%%%%%%%%%%%%%%%%%%%%%%%%%%%%%%%%%%%%%%%%%%%%%%%%%%%%%%%%%%%%%%%%%%%%%%%%%%%%%%%%%%%%
\section{Motivating Example}
\label{sec:motivation}
% %%%%%%%%%%%%%%%%%%%%%%%%%%%%%%%%%%%%%%%%%%%%%%%%%%%%%%%%%%%%%%%%%%%%%%%%%%%%%%%%%%%%%%%%%%%%%%%%%%%%%%%%%%%%%
% %%%%%%%%%%%%%%%%%%%%%%%%%%%%%%%%%%%%%%%%%%%%%%%%%%%%%%%%%%%%%%%%%%%%%%%%%%%%%%%%%%%%%%%%%%%%%%%%%%%%%%%%%%%%%
One of the most important scientific discoveries of the past century was undoubtedly the observation of the cosmic microwave background (CMB); see \cite{penzias64}. The discovery is one of the strongest pieces of evidence for the Big Bang theory, which predicted its existence.  The properties of the CMB carry important information about how the universe has evolved from very early epochs, and therefore there is great interest in making a full sky measurement of CMB as accurately as possible. To date, three satellites --- COBE,  WMAP  and most recently Planck --- have done this with increasing sensitivity and resolution. Planning for a future mission that would make even more detailed observations has started \cite{delabrouille18}.

An important problem is that the signals measured by these satellites do not contain only CMB radiation but also contributions from a number of other microwave sources. Hence the contribution of CMB to the observed microwave brightness at any point in the sky must be separated out. These other sources come in two types.  There are point sources, consisting of stars and galaxies; these are relatively easy to identify and remove from the data.  There are also diffuse sources, such as microwaves emitted from dust in our galaxy.  These are harder to separate from the CMB, which is also a diffuse signal. The source separation problem concerns these diffuse sources only, assuming that the point sources have been removed prior to analysis.

% %%%%%%%%%%%%%%%%%%%%%%%%%%%%%%%%%%%%%%%%%%%%%%%%%%%%%%%%%%%%%%%%%%%%%%%%%%%%%%%%%%%%%%%%%%%%%%%%%%%%%%%%%%%%%
\subsection{Source separation model for all-sky maps}
This source separation problem can be formulated as follows. There are $n$ all-sky maps of microwave brightnesses that have been observed at frequencies $\nu_1,\ldots,\nu_n$.  These maps consist of a discrete grid of values over $N$ pixels that partition the entire sky. In astronomy, the HEALPix partition into pixels of equal area is used \cite{healpix05}.  HEALPix partitions the sky into 12 'base' areas, and a separate analysis is done on each; see Appendix \ref{app:healpix} for more details. For Planck data, $n=9$ and there are 12 base areas each with $N = 4^{10} = 1,048,576$ pixels, for a total of 12,582,912 pixels across the entire sky.  The following analysis is repeated 12 times for each base patch and the results are stitched together.

Let $\vek y_j \in \mathbb{R}^n$ denote the brightness of the maps at pixel $j$. Assuming that point sources have been dealt with, the brightness at each point in the sky is a combination of the $m$ different diffuse sources.  Letting $\vek s_j \in \mathbb{R}^m$ denote the vector of source brightnesses at pixel $j$, the source separation model takes the form:
\begin{equation} 
\vek y_j = \vek A \vek s_j + \beps_j, \; j=1,\ldots,N,
\label{eq:model}
\end{equation}
where $\vek A$ is an $n \times m$ matrix of weights (known as the mixing matrix in source separation) and $\beps_j$ is a vector of independent measurement errors, assumed to be independently Gaussian distributed with mean 0 and known variances $\sigma_1^2,\ldots,\sigma_n^2$.  
% The CMB is assumed to be the first component of $s_j$. The components of the mixing matrix $A_{kl}$ can be interpreted as the relative contribution of the $l$th source at the $k$th frequency. 
 In the full inference task, $\vek A$ is defined by a small number of parameters whose values must be learned from the data. Because that aspect of the learning task is not relevant to this work, here the elements of $\vek A$ are assumed to have known values. Appendix \ref{app:A} gives the details of the construction of $\vek A$.

% This model also specifies the factor analysis and independent components analysis models from statistics and machine learning {\tt <<citation>>}.  In the former case $A$ is known as the factor loading matrix.  

% {\color{red} {\tt Comment KMS: it might be good to clarify the shapes of these matrices we are constructing in the next sentences for clarity. DONE --- SPW
% }}
Let $\bcalY_k = \left(y_{1k},\ldots,y_{Nk}\right)$ denote the $k$th map (in other words, the $k$th component of each $\vek y_j$) and $\vek Y = \left(\bcalY_1,\ldots,\bcalY_n\right) \in \R^{nN}$ denote the vector of all observations, stacked by map. Similarly $\bcalS_l = (\vek s_{1l},\ldots,\vek s_{Nl})$ is the map of the $l$th source and $\vek S = (\bcalS_1,\ldots,\bcalS_m) \in \R^{mN}$. Following Equation \ref{eq:model}, one has
\begin{equation*} 
\vek Y = \vek B \vek S + \bcalE,
% \label{eq:model_stacked}
\end{equation*}
where $\vek B = \vek A \otimes \vek I_N$ is an $nN \times mN$ matrix, $\vek I_N$ is the $N$ dimensional identity matrix, and $\bcalE$ is the vector of $nN$ Gaussian distributed measurement errors with mean vector 0 and $nN \times nN$ dimensional diagonal precision matrix $\vek C$, with the first $N$ entries on the diagonal equal to $1/\sigma_1^2$, the next $N$ equal to $1/\sigma_2^2$, and so on.  Hence, given $\vek A$, $\vek S$ and the $\sigma_k^2$, $\vek Y$ is Gaussian distributed with mean $\vek B \vek S$ and precision matrix $\vek C$:
\begin{equation}
    p(\vek Y \, | \, \vek S) = (2\pi)^{Nn/2} |\vek C|^{1/2} \, \exp \left\{ - \frac{1}{2} (\vek Y - \vek B \vek S)^T \vek C (\vek Y - \vek B \vek S) \right\}.
\label{eq:likelihood}
\end{equation}
Note that, for Planck data,  $nN = 9 \times 4^{10} = 9,437,184$ 
% 113,246,208
and typically one uses $m=4$ diffuse sources so that $mN = 4 \times 4^{10} = 4,194,304$. % 50,331,648.

The final element of the problem is a model for the diffuse sources.  In a Bayesian learning approach, this model is a probability distribution $p(\vek S)$, known as the prior distribution, that serves to regularise the separation task. In our approach, this distribution models the spatial smoothness of each diffuse source through the use of an intrinsic Gaussian Markov random field model \cite{rue05}. This model ensures that the intensity of a source at neighbouring pixels in the sky cannot be too different; more concretely, if pixel $j$ and $j^{\prime}$ are neighbours then $y_{jk} - y_{j^{\prime}k}$ is zero-mean Gaussian distributed with a variance $\varphi_k$ whose value must be specified. The different source maps are assumed to be independent. This leads to a zero-mean multivariate Gaussian form for $p(\vek S)$ with a block diagonal form for the covariance matrix $\vek Q^{-1}$, one block being associated with each source:
\begin{equation}
    p(\vek S) \: = \: (2\pi)^{Nm/2} |\vek Q|^{1/2} \, \exp \left\{ - \frac{1}{2} (\vek Y - \vek B \vek S)^T \vek Q (\vek Y - \vek B \vek S) \right\},
    \label{eq:prior}
\end{equation}
where $\vek Q = \mbox{bldiag}(\vek{Q}_1(\varphi_1),\ldots,\vek{Q}_m(\varphi_m))$ is an $mN \times mN$ matrix and the $\vek{Q}_l(\varphi_l)$ are sparse, positive semi-definite $N \times N$ matrices, each depending on a single parameter $\varphi_l$. Appendix \ref{app:igmrf} has details on the form and derivation of the $\vek{Q}_l$.

% %%%%%%%%%%%%%%%%%%%%%%%%%%%%%%%%%%%%%%%%%%%%%%%%%%%%%%%%%%%%%%%%%%%%%%%%%%%%%%%%%%%%%%%%%%%%%%%%%%%%%%%%%%%%%
\subsection{Source separation solution}
The primary task is to use the data $\vek Y$ to learn about all sources $\vek S$.  A Bayesian learning approach is adopted, where the goal is to derive the so-called posterior probability distribution $p(\vek S \, | \, \vek Y)$.  Following Bayes' Law:
\[ p(\vek S \, | \, \vek Y) \: \propto \: p(\vek Y \, | \, \vek S) \: p(\vek S), \]
where Equations \ref{eq:likelihood} and \ref{eq:prior} specify $p(\vek Y \, | \, \vek S)$ and $p(\vek S)$. 

% Tackling the computation of $p(\vek S, \psi \, | \, \vek Y)$ gives rise to the linear system. By the multiplication law of probability, we have:
% \[ p(\vek S, \psi \, | \, \vek Y) \: = \: p(\psi \, | \, \vek Y) \: p(\vek S \, | \,  \psi, \vek Y), \]
% We show in {\tt <<citation>>} that if one can compute the term $p(\vek S \, | \,  \psi, \vek Y)$ then there is a feasible approach to computing $p(\psi \, | \, \vek Y)$.  This also allows Bayesian estimates of the model parameters and sources, including CMB, to be derived. A big advantage of the approach is that the distribution allows one to quantify uncertainty in the estimate, an important factor in understanding how much one can reasonably infer from the data. 
Some simple properties of the Gaussian distribution show that $p(\vek S \, | \,  \vek Y)$ is another Gaussian with precision matrix
\begin{equation}
   \widehat{\vek Q} = \vek B^T \vek C \vek B + \vek Q
    \label{eq:Q*}
\end{equation}
and mean vector $\widehat{\bmu}$ that satisfies
\begin{equation}
\widehat{\vek Q}\, \widehat{\bmu} = \vek B^T \vek C \vek Y.
\label{eq:mu*}
\end{equation}
So key to implementing the source separation is the solution to Equation \ref{eq:mu*}, under the structures for $\vek Q$, $\vek B$ and $\vek C$ that have been described, as this gives a point estimate of the sources.  In addition, computing elements of $\widehat{\vek Q}$ (of most importance being the main diagonal that contains the variance of the sources at each pixel) gives information about uncertainty of and correlation between source values across the sky.  We elaborate on implementing the action of $\widehat{\vek Q}$ in the subsequent section.

% \[ p(\psi \, | \, \vek Y) \: \propto \:  \frac{p(\vek Y \, | \, \psi, \vek Y) \: p(\vek S \, | \, \psi) \: p(\psi)}{p({\cal S \, | \, \psi,\\vek Y)} \]

% %%%%%%%%%%%%%%%%%%%%%%%%%%%%%%%%%%%%%%%%%%%%%%%%%%%%%%%%%%%%%%%%%%%%%%%%%%%%%%%%%%%%%%%%%%%%%%%%%%%%%%%%%%%%%
% %%%%%%%%%%%%%%%%%%%%%%%%%%%%%%%%%%%%%%%%%%%%%%%%%%%%%%%%%%%%%%%%%%%%%%%%%%%%%%%%%%%%%%%%%%%%%%%%%%%%%%%%%%%%%
\section{Solving the Linear System}
\label{sec:solution}
% %%%%%%%%%%%%%%%%%%%%%%%%%%%%%%%%%%%%%%%%%%%%%%%%%%%%%%%%%%%%%%%%%%%%%%%%%%%%%%%%%%%%%%%%%%%%%%%%%%%%%%%%%%%%%
% %%%%%%%%%%%%%%%%%%%%%%%%%%%%%%%%%%%%%%%%%%%%%%%%%%%%%%%%%%%%%%%%%%%%%%%%%%%%%%%%%%%%%%%%%%%%%%%%%%%%%%%%%%%%%

%The task is to solve Equation \ref{eq:mu*} for $\mu^*$ where $Q^* = B^T C B + Q$ and:
%\begin{itemize}
%    \item $B$ is an $nN \times mN$ matrix of the form $A \otimes I_N$, where $A$ is a small matrix ($n \times m$) and $N$ can be large (tens of millions in the case of our application). Therefore it is sparse but the non-zero elements are located across the matrix and not in a band near the diagonal.
%    \item $C$ is an $nN \times nN$ diagonal matrix with strictly positive elements on the diagonal.
%    \item $Q$ is an $mN \times mN$ block-diagonal matrix, with $m$ blocks, each block being a sparse, semi-positive definite matrix of dimension $N \times N$.
%\end{itemize}

Although $\vek Q$ is only semi-positive definite, $\widehat{\vek Q}$ will be positive definite.  As is often the case, the sparsity of $\widehat{\vek Q}$ does not guarantee that its Cholesky decomposition will be sparse. %; indeed, we can demonstrate this for a small sample realization; see Figure \ref{fig:sparse}.  
One could consider sparse direct methods (see, e.g., \cite{GilbertSchreiber:1992:Sparse-parallel-cholesky}) in this setting, but this problem has so much additional structure that we explore matrix-free approaches which take advantage of Kronecker structure.  Thus, we do not consider sparse-direct methods in this paper, though we do use SuperLU \cite{superlu_smp99, superlu_ug99, superlu-userguide-li05} as an inner solver for a method suggested in \cite{kohller} against which we test our proposed solvers, cf. \Cref{section:comparison-Kohller}.  We also note that our descriptions of quantities arising in Kronecker form and the use reshaping to map between Kronecker form and unvectorised form is \emph{only a matter of presentation}.  We do not advocate reshaping to go between Kronecker and unvectorised forms.

% \begin{figure}
% \centering
% \includegraphics[width=50mm]{Qsparse.jpg} 
% \includegraphics[width = 50mm]{LofQ.jpg}
% \caption{\label{fig:sparse}Sparseness of $\widehat{\vek Q} = \vek Q + \vek B^T \vek C \vek B$ (left) and it Cholesky decomposition (right) for a small example ($m=4$ and $N=256$).  Non-zero elements are shaded.}
% \end{figure}

The problem being solved is of the form
\be\label{eqn.Axb}
	\prn{\vek Q + \vek B^{T}\vek C\vek B}\boldsymbol\mu = \vek B^{T}\vek C\vek Y\mwith \prn{\vek Q + \vek B^{T}\vek C\vek B}\in\R^{mN\times mN}
\ee
where $m\in\curl{4,5}$ and $N=\CO(10^{7})$. As described in Appendix \ref{app:igmrf}, $\vek Q$ is a $m$-block diagonal matrix with block $i$ being of the form $\vek{Q}_i(\varphi_i) = \varphi_{i}\vek D^{T}\vek D=\varphi_{i}\vek D^{2}$, for $\varphi_i > 0$.  
The matrix $\vek D\in\R^{N\times N}$ is symmetric and its elements encode a nearest neighbor 
coupling of nodes from the discretization of the sky, whereby
\be\nn
	\vek D = \prn{D_{ij}} \mwith D_{j_{1},j_{2}}=\begin{cases} 1 & \mbox{nodes $j_{1}\neq j_{2}$ are neighbors}\\ 
																								0 & \mbox{nodes $j_{1}\neq j_{2}$ are not neighbors}\\
																								-\sum_{\underset{\ell\neq j}\ell=1}^{N}D_{j\ell}& j:=j_{1}=  j_{2}\end{cases}.
\ee
The matrix $\vek B\in\R^{nN\times mN}$ with $n=9$ can be written as the Kronecker product $\vek B = \vek A \otimes \vek I_{N}$
where $\vek A\in\R^{n\times m}$. We also have that $\vek C\in\R^{nN\times nN}$ is a diagonal matrix with all positive entries.  The matrix $\vek C$ also has Kronecker structure; namely, there 
are $N$ positive integer values $n_1,\ldots,n_N$ and $n$ positive values $\tau_1,\ldots,\tau_n$ 
\be\nn
	\vek C = {\rm diag}\curl{\tau_{1},\tau_{2},\ldots,\tau_{n}} \otimes {\rm diag}\curl{n_{1},n_{2},\ldots, n_{N}} =: \vek T\otimes \vek N.
\ee

%Because of the large scale of this problem and sparse, Kronecker structure, \cref{eqn.Axb} should not be solved directly by a matrix decomposition method.
%Instead 
We propose treating this problem using iterative methods built on efficient implementations of matrix-vector products (matvecs) for the matrices
from this problem (cf. \cref{section.CG-strategy}) or as the low-rank approximation to the solution of a 
Sylvester matrix equation (cf. \cref{section.Sylvester-strategy}).  Furthermore, most of the system matrix should never be explicitly constructed.  Rather, it should only be represented as a stored procedure
implementing matrix-vector product.  The only matrix explicitly represented should be $\vek A$.

% %%%%%%%%%%%%%%%%%%%%%%%%%%%%%%%%%%%%%%%%%%%%%%%%%%%%%%%%%%%%%%%%%%%%%%%%%%%%%%%%%%%%%%%%%%%%%%%%%%%%%%%%%%%%%
\subsection{Efficient matrix-vector products}
\label{subsec:matvec}
We describe first how to perform the matrix-vector product for each part of the full system matrix and then we finish by combining them, yielding a routine
for $\vek x \mapsto \prn{\vek Q + \vek B^{T}\vek C\vek B}\vek x$.  This will be necessary for using an appropriate iterative method whose core operation
is this matvec routine. 
% {\color{red}[KMS: need to remove explicit Kronecker product stuff not used in code that confused referees]}
% \subsubsection{Matvec routine for $\vek v\mapsto\vek D\vek v$}
Observe that since $\vek D$ is a sparse matrix encoding the nearest-neighbor coupling for each node, the routine $\vek v \mapsto\vek D\vek v$ involves
 updating entry $j$ of $\vek v$ (associated to node $j$ of grid) using its four nearest neighbors and the value at the note itself. Thus a sequential matrix-free routine can be easily formulated.
 \RestyleAlgo{boxruled}
\begin{algorithm2e}
\caption{Matvec procedure $\vek v\mapsto \vek D\vek v$ \label{alg.MatvecD}}
Given: vector $\vek v\in\R^{N}$, perhaps represented as nodal values on the grid \\
\For{$j=1,2,\ldots N$}{
	$\vek v(j) \leftarrow \vek v(j_{right}) + \vek v(j_{left}) + \vek v(j_{above}) + \vek v(j_{below}) - \prn{\sum_{\underset{\ell\neq j}\ell=1}^{N}D_{j\ell}}\vek v(j)$
}
\end{algorithm2e}
For Algorithm \ref{alg.MatvecD}, we note that since the matrix $\vek D$ does not change, we can pre-compute and store the rowsums $\sum_{\underset{\ell\neq j}\ell=1}^{N}D_{j\ell}$.  Also, it should be noted that Algorithm \ref{alg.MatvecD} is parallelizable, as the grid can be
divided up and sent to different processes, and communication between processors only needs to occur for neighboring nodes that are on different processors. 
The computational cost of a serial version of Algorithm \ref{alg.MatvecD} is $\CO(N)$, since we do $4$ floating point operations per entry of $\vek v$.

% \subsubsection{Matvec routine for $\vek u\mapsto\vek Q\vek u$}
% \label{subsec:matvecQ}
As already discussed, $\vek Q$ is a block diagonal matrix with blocks $\varphi_{i}\vek D^{2}$.  Thus, it can be represented as a Kronecker product, namely
\be\nn
	\vek Q = \begin{bmatrix}
	\varphi_{1} &  &  &  \\ 
	 & \varphi_{2} &  &  \\ 
	 &  & \ddots &  \\ 
	 &  &  & \varphi_{n}
	\end{bmatrix} \otimes \vek D^{2}=: \vek P\otimes\vek D^{2}.
\ee
For $\vek u = \vek w\otimes \vek v\in\R^{nN}$ with $\vek w\in\R^{n}$, we  
perform this matrix-vector product by instead considering a reshaping of $\vek u$.  
Using the \matlab{-style} \texttt{reshape} notation, let $\vek U = {\tt reshape}(\vek u, N, n)\in\R^{N\times n}$; i.e., 
	\be\nn
		\mif\vek u = \bbmat \vek u_{1}\\ \vek u_{2}\\\vdots \\ \vek u_{n} \ebmat \mthen \vek U = {\tt reshape}(\vek u, N, n) = \bbmat \vek u_{1} & \vek u_{2} & \cdots & \vek u_{n} \ebmat.
	\ee

Then the matvec in  Algorithm \ref{alg.MatvecD} can 
also be represented by reshaping $\vek u$ and computing $\vek D^{2}\vek U\vek P$.  
% \subsubsection{Matvec routine for $\vek x\mapsto \vek B^{T}\vek C\vek B\vek x$}
Again, we take advantage of the Kronecker structures of $\vek B$ and $\vek C$.   We can write
  \be\label{eqn.MatvecBCBT-kron}
  	\vek B^{T}\vek C\vek B\vek x = \prn{\vek A^{T}\otimes \vek I_{N}}\prn{\vek T\otimes \vek N}\prn{\vek A\otimes \vek I_{N}}\vek x.
  \ee  
  If we set $\vek M =  {\tt reshape}(\vek x, N, n)$, then the matrix-vector product can be efficiently executed via the block vector mapping
  \begin{align}
      \vek M
      \mapsto 
      \prn{\vek N\prn{\vek M\vek A^T}\vek T}\vek A,
      \label{eq:block-vector-mapping-matvec}
  \end{align}
  which enables the running of a matrix-free Krylov iteration on a Kronecker-structured system implicitly, without forming any Kronecker products.
\subsection{Solution of \cref{eqn.Axb} using matrix-free iterative methods}\label{section.CG-strategy}
The first strategy we propose entails considering the linear system \cref{eqn.Axb} as it is presented with the 
coefficient matrix being built from and applied implicitly using the Kronecker product formulation.  We propose to use
a so-called ``matrix-free'' method which only requires a stored procedure to perform
$\vek x \mapsto \prn{\vek Q + \vek B^{T}\vek C\vek B}\vek x$, implemented in block form via \eqref{eq:block-vector-mapping-matvec}, meaning we need no explicit representation of
the coefficient matrix.  

We briefly review standard Krylov subspace methods and introduce some terminology.
For a more detailed description, see, e.g.,
\cite{Saad.Iter.Meth.Sparse.2003,szyld.simoncini.survey.2007,vdVorst} and the references therein.
Consider solving
\begin{equation}\label{eqn.correction}
	\vek G(\vek x_{0} + \vek t) = \vek b
\end{equation}
where $\vek G\in\Rnn$, $\vek x_{0}\in\Rn$ is an initial approximation, and $\vek t$ is the unknown correction.
Let $\vek r_{0} = \vek b - \vek G\vek x_{0}$ be the initial residual.
In iteration $j$, a Krylov method selects a correction $\vek t_{j}$ from the
$j$th Krylov subspace generated by $\vek G$ and
$\vek r_{0}$,
\begin{equation}\label{eqn.kryl-def}
	\vek t_{j} \in \CK_{j}(\vek G,\vek r_{0}) = {\rm span}\left\lbrace \vek r_{0}, \vek G\vek r_{0}, \vek G^{2}\vek r_{0}\cdots, \vek G^{j-1}\vek r_{0} \right\rbrace
	=\left\lbrace p(\vek G)\vek r_{0}\ |\deg p < j\right\rbrace,
\end{equation}
i.e., from the space of polynomials of degree less than $j$ in $\vek G$ acting on $\vek r_{0}$.
The constraint space defines
the Krylov subspace method up to implementation choices and must also be compatible with the correction space
such that the method being implicitly defined can be efficiently implemented and exhibits good convergence characteristics.

For coefficient matrices $\vek G$ that are
symmetric positive definite (such as that in \cref{eqn.Axb}), 
the celebrated, efficient, and stable method of Conjugate Gradients is the 
Krylov subspace method of choice.  Implicitly, at iteration $j$, the method has generated the subspace 
$\CK_{j}\prn{\vek G,\vek r_{0}}$ and computes $\vek t_{j}\in\CK_{j}\prn{\vek G,\vek r_{0}}$ such that
the approximation $\vek x_{j} = \vek x_{0} + \vek t_{j}$ minimizes 
\be\nn
    \nm{\vek x - \vek x_j}^2_{G}
    =
	\nm{\vek b - \vek G\vek x_{j} }^2_{G^{-1}} 
    = 
    \prn{\vek b - \vek G\vek x_{j} }^{T}\vek G^{-1}\prn{\vek b - \vek G\vek x_{j} }.
\ee
The actual derivation of the method can be found in, e.g., \cite{Saad.Iter.Meth.Sparse.2003}.  
We note that we have omitted a Lanczos-based derivation of CG, as its derivation is well documented in the literature.  We derive an alternative approach based on block Lanczos in \Cref{section.Sylvester-strategy}.
We present CG for the case
$\vek G =  \prn{\vek Q + \vek B\vek C\vek B^{T}}$ and right-hand side 
$\vek b = \vek B^T\vek C\vek Y$ as Algorithm \ref{alg.CGnoprecond}.
%\CGnoprecond{1111}{2222}{3333}{4444}{55555}
\begin{algorithm2e}
\caption{Method of Conjugate Gradients -- for $\vek G =  \prn{\vek Q + \vek B\vek C\vek B^{T}}$ and $\vek b = \vek B^T\vek C\vek Y$\label{alg.CGnoprecond}}
Given: $\vek y$, $\vek B$, $\vek C$, and $\vek Q$ defined as above and $\vek x_{0}$; $\varepsilon>0$ a convergence tolerance\\
$\vek x = \vek x_{0}$; $\vek r = \vek B\vek C\vek y - \prn{\vek Q + \vek B\vek C\vek B^{T}}\vek x$; $\vek p = \vek r_{0}$\\
\While{$\nm{\vek r}  >\varepsilon\nm{\vek r_{0}}$}{
	$\vek r_{old}\leftarrow \vek r$\\
	$\vek q\leftarrow \prn{\vek Q + \vek B\vek C\vek B^{T}}\vek p$\\
	$\alpha = \frac{\vek r^{T}\vek r}{\vek p^{T}\vek q}$\\
	$\vek x \leftarrow \vek x + \alpha\vek p$\\
	$\vek r\leftarrow \vek r - \alpha\vek q$\\
	$\beta \leftarrow \frac{\vek r^{T}\vek r}{\vek r_{old}^{T}\vek r_{old}}$\\
	$\vek p \leftarrow \vek r + \beta\vek p$
}
\end{algorithm2e}

\textbf{Advantages.} The CG methods has a number of advantages in this setting.
The coefficient matrix of \cref{eqn.Axb} does not need to be constructed (and there is no expensive decomposition).
The convergence speed is determined by the eigenvalues of the coefficient matrix and the specific right-hand side.  If convergence speed is too slow, \emph{preconditioning} techniques can be used to accelerate convergence.  Furthermore, in this setting CG is quite amenable to parallelization by distributing sky domain grid onto different processors and parallelizing
	matvec procedures.

\textbf{Disadvantages.} The CG method also has some relevant disadvantages in this setting. We still must deal with vectors of size $nN$. We also cannot take advantage of any structure or compressibility which may be present in $\vek Y$. Most importantly, it has been shown that applying an iterative method directly to the tensor-structured problem can lead to performance hampered by the large condition number of the full-size problem \cite{pal-sim-kroncond}.

Understanding both the advantages and disadvantages, we have been able to demonstrate for this problem using real data that the CG method can be used effectively to solve the problem.  One observation we make is that the coefficient matrix arises from the Bayesian regularisation of the source separation problem.  Such matrices are often more well-conditioned due to the regularisation.
% %%%%%%%%%%%%%%%%%%%%%%%%%%%%%%%%%%%%%%%%%%%%%%%%%%%%%%%%%%%%%%%%%%%%%%%%%%%%%%%%%%%%%%%%%%%%%%%%%%%%%%%%%%%%%
\subsection{Solution of \cref{eqn.Axb} using a Sylvester interpretation}\label{section.Sylvester-strategy}
Going one step further, we can simply reshape \cref{eqn.Axb} and interpret directly as a Sylvester 
matrix equation, which will allow us to take advantage of additional structure in $\vek y$ that may be available.

% {\color{red} SPW comment: in general $\vek y$ can be any set of values in $\R^{nN}$ (although for the CMB application we'd usually have $y_i > 0$). So it does not in general have any kind of Kronecker structure for example. It is the concatenation of $n$ images each of size $N$, which reshaping to an $N \times n$ matrix emphasises.  Not clear if we need to assume further structure on $\vek y$ for the Sylvester approach?}

We can reshape the solution ($\vek M =  {\tt reshape}(\boldsymbol\mu, N, n)$) and ($\vek Y =  {\tt reshape}(\vek y, N, m)$).
With these reshapings, we can rewrite \cref{eqn.Axb} as the generalize Sylvester equation 
\be\nn
	\vek D^{2}\vek M\vek P + \vek N\vek M\vek A^{T}\vek T\vek A = \vek N\vek Y\vek T\vek A.
\ee
As $\vek N$ and $\vek P$ are diagonal and invertible, we can easily multiply on the left by $\vek N^{-1}$ and on the
right by $\vek P^{-1}$, yielding the Sylvester equation
\be\label{eqn.Axb-Sylv}
	\vek N^{-1}\vek D^{2}\vek M + \vek M\vek A^{T}\vek T\vek A\vek P^{-1} = \vek Y\vek T\vek A\vek P^{-1}.
\ee
A survey by Simoncini covering a variety of solution techniques was recently published \cite{Simoncini.MatEq-Survey.2016}.
Since $\vek A^{T}\vek T\vek A\vek P^{-1}\in\Rnn$, this falls into a category of Sylvester equations with one 
large, sparse coefficient matrix ($\vek N^{-1}\vek D^{2}$) and one small, dense one ($\vek A^{T}\vek T\vek A\vek P^{-1} $). Following \cite{kohller},  we denote this as \emph{sparse-dense} category. We describe a technique from the literature wherein this structure has been previously exploited.

\subsubsection{An existing technique for the treating sparse-dense Sylvester equations}
Since \eqref{eqn.Axb-Sylv} falls into this category, it is of the form
\begin{equation}
\label{eq:sparse-dense-Sylvester}
    \mathbf{HM} + \mathbf{ML} = \mathbf{F},
\end{equation}
with $\vek H$ large and sparse and $\vek L$ small and dense.
In \cite{kohller}, the authors propose Algorithm \ref{alg:big-small-Kohler} to solve \eqref{eq:sparse-dense-Sylvester}.  In the numerical experiments (cf. \Cref{section:comparison-Kohller}), we demonstrate that this algorithm is not appropriate for the full-scale application problem we are treating in this manuscript.  For small problems, it does produce solutions with smaller residual norms than what we propose in the subsequent section, but it exhibits quadratic growth in solve time as we increase the problem size, mainly due to the cost of inner linear solves required by method; see \Cref{section:comparison-Kohller} for further details.  

We propose a different approach, which we demonstrate (for test problems involving real Planck data) is more robust with respect to the increase in run time as the size of the problem (i.e., size of the data) increases, and we demonstrate that it delivers solutions with sufficient accuracy in a timely manner.  It uses the fact that equations with Sylvester structure can be treated using block Krylov subspaces, exploiting a block generalization of the shift-invariance exhibited by Krylov subspaces; see, e.g., \cite{Freund:1993:1} concerning shift invariance and \cite{Simoncini.MatEq-Survey.2016} and references therein for its block generalization to the Sylvester equation setting.

\begin{algorithm2e}
    \caption{Solution of the Sparse-Dense Sylvester Equation (from \cite{kohller}) \label{alg:big-small-Kohler}}
    \textbf{Input:} $\mathbf{H} \in  \mathbb{R}^{N \times N}, \vek L \in \mathbb{R}^{n \times n}$ and $\vek F \in \mathbb{R}^{\mathbf{N} \times n}$ defining \eqref{eq:sparse-dense-Sylvester}\\
    \textbf{Output: $\mathbf{M} \in \mathbb{R}^{N \times n}$ solving \eqref{eq:sparse-dense-Sylvester}}\\
    Compute Schur decomposition $\vek L = \vek W\vek R\vek W^{T}, \vek W^{T}\vek W = \vek I_{n}$ \\
  $\widehat{\vek F} := {\vek F\vek W}$ \\
  \For{$i=1,2,...,n$}{
    $\widehat{\vek F} := {\vek F}_{i} - \sum^{i-1}_{k=1}\widehat{\vek X}_{k}R_{ki}$ \\
     Solve $(\vek H + R_{ii}\vek I)\widehat{\vek X}_{i} = \widehat{\vek F}$\\
  }
  $\vek M := \widehat{\vek X}\vek W$ \\
  \textbf{return} $\vek M$
\end{algorithm2e}

\subsubsection{Application of block Krylov approach for large-scale problems}

Large-scale Sylvester-structured problems of sparse-dense type are amenable to solution by 
generating a \emph{block} Krylov subspace associated to the large coefficient matrix, i.e., $\vek N^{-1}\vek D^{2}$.
This space at iteration $j$ has the form
\bea
	\K_{j}\prn{\vek N^{-1}\vek D^{2}, \vek Y\vek T\vek A\vek P^{-1}} &:=& \CK_{j}\prn{\vek N^{-1}\vek D^{2}, \vek Y\vek T\vek A\vek e_{1}} \nn\\&&+ \CK_{j}\prn{\vek N^{-1}\vek D^{2}, \vek Y\vek T\vek A\vek e_{2}} + \cdots\nn\\ &&+ \CK_{j}\prn{\vek N^{-1}\vek D^{2}, \vek Y\vek T\vek A\vek e_{m}},\nn
\eea
where $\vek e_{i}$ is the $i$th Cartesian basis vector.  The Sylvester equation is projected onto this space, and a smaller Sylvester problem is solved at each iteration.  The matrix 
$\vek N^{-1}\vek D^{2}$ is no longer symmetric positive definite in the standard inner product; however, it is if we consider 
a different inner product. 

Any symmetric positive-definite matrix can be used to define a new inner product.  The matrix $\vek N$ is a diagonal matrix with 
only positive entries on its diagonal; therefore, it is a positive-definite matrix and can be used to induce the inner product $\prn{\vek x,\vek y}_{\vek N}= \vek y^{T}\vek N\vek x$.  A matrix $\vek G$ is symmetric with respect to $\prn{\cdot,\cdot}_{\vek N}$ if for all $\vek x,\,\vek y\in\Rn$
$\prn{\vek G\vek x,\vek y}_{\vek N} = \prn{\vek x,\vek G\vek y}_{\vek N}$; i.e., $\vek y^{T}\vek N\vek G\vek x= \vek y^{T}\vek G^{T}\vek N\vek x$.   

What is important here is that we exploit the property proven in the following Lemma. 
\begin{lemma}
The coefficient matrix $\vek N^{-1}\vek D^{2}$ is symmetric with respect to $\prn{\cdot, \cdot}_{\vek N}$.
\end{lemma}
\bproof
	This is straightforward to show by direct computation.  For all $\vek x,\, \vek y\in\Rn$, we have 
	\be\nn
		\prn{\vek N^{-1}\vek D^{2}\vek x,\vek y}_{\vek N} = \vek y^{T}\vek N\vek N^{-1}\vek D^{2}\vek x = \vek y^{T}\vek D^{2}\vek x = \vek y^{T}\vek D^{2}\vek N^{-1}\vek N\vek x=\prn{\vek x,\vek N^{-1}\vek D^{2}\vek y}_{\vek N}.
	\ee  
	Thus the symmetry with respect to $\prn{\cdot, \cdot}_{\vek N}$ is proven.
\eproof

Iterative methods designed for symmetric matrices can be
converted to methods for matrices symmetric with respect to any inner product simply by using that inner product in the algorithm.
We can implement a Conjugate Gradient-like iteration which uses
a fixed amount of memory, is amenable to parallelization, and can be implemented to have good 
data movement characteristics \footnote{For problems of this size, we must concern ourselves with how much 
data an algorithm moves, as this becomes a dominating cost.}
Thus it behooves us to elaborate on block Lanczos-based methods for symmetric matrices.

% %%%%%%%%%%%%%%%%%%%%%%%%%%%%%%%%%%%%%%%%%%%%%%%%%%%%%%%%%%%%%%%%%%%%%%%%%%%%%%%%%%%%%%%%%%%%%%%%%%%%%%%%%%%%%
\subsubsection{The block Lanczos process with respect to $\prn{\cdot, \cdot}_{\vek N}$} We mentioned in \Cref{section.CG-strategy} that 
CG is a Krylov subspace method.  However, we avoided describing any of the details about how this 
subspace is used to derive CG for brevity.  We elaborate here on some aspects of Krylov subspaces for symmetric
coefficient matrices to derive an iterative method for Sylvester equations whose larger coefficient matrix (here $\vek N^{-1}\vek D^{2}$) is symmetric with respect to respect to $\prn{\cdot, \cdot}_{\vek N}$.

Let us denote $\widehat{\vek Y} = \vek Y\vek T\vek A\vek P^{-1}$ and $\widehat{\vek L} = \vek A^{T}\vek T\vek A\vek P^{-1} $.
 The symmetric Lanczos process is an iterative procedure to generate an $\vek N$-orthonormal basis for $\K_{j}\prn{\vek N^{-1}\vek D^{2}, \widehat{\vek Y}} $ block-by-block.  
 What has been shown is that if the matrix is
 symmetric with respect to the inner product used for orthogonalization, at step $i$, the $(i+1)$st block of vectors generated at that step are already orthogonal with respect to all but blocks $i$ and $i-1$, meaning 
 the orthogonalization procedure requires we only store the two most recently generated blocks.  At step $i$, we have generated the orthonormal basis $\curl{\vek V_{1},\vek V_{2},\cdots, \vek V_{i}}$ with 
 $\vek V_{\ell}\in\R^{n\times m}$, $\vek V_{\ell}^{T}\vek N\vek V_{k}=\vek 0$ for $\ell\neq k$, and $\vek V_{\ell}^{T}\vek N\vek V_{\ell}=\vek I_{m}$, the $m\times m$ identity matrix.  To obtain the next basis vector, 
 we compute $\vek W_{i+1} = \prn{\vek N^{-1}\vek D^{2}}\vek V_{i}$.  Normally, we would then need to orthogonalize these vectors against all previous, but since $\vek N^{-1}\vek D^{2}$ is symmetric with respect to
 $\prn{\cdot,\cdot}_{\vek N}$, it has been proven that $\vek V_{\ell}^{T}\vek N\vek W_{i+1}=\vek 0$ for $\ell < i-1$.  The orthogonalization coefficients are stored in the matrix $\underline{\vek T}_{j}\in\R^{(j+1)m,jm}$ 
 whose first $jm$ rows are the symmetric, block tridiagonal matrix $\vek T_{j}\in\R^{jm\times jm}$.
These have the structure
 \be\nn
 	\underline{\vek T}_{j} =  	
 	\bbmat 
	 	\vek H_{1} & \vek B_{2} &  &  &  \\ 
 		\vek B_{2} & \vek H_{2} & \vek B_{3} &  &  \\ 
 		 & \vek B_{3} & \ddots & \ddots &  \\ 
 		 &  & \ddots & \vek H_{j-1} & \vek B_{j} \\ 
 		 &  &  & \vek B_{j} & \vek H_{j} \\ 
 		 &  &  &  & \vek B_{j+1} 	
 	\ebmat \mand 
 	\vek T_{j} =  	
 	\bbmat 
	 \vek H_{1} & \vek B_{2} &  &  &  \\ 
 		\vek B_{2} & \vek H_{2} & \vek B_{3} &  &  \\ 
 		 & \vek B_{3} & \ddots & \ddots &  \\ 
 		 &  & \ddots & \vek H_{j-1} & \vek B_{j} \\ 
 		 &  &  & \vek B_{j} & \vek H_{j} 
 	\ebmat,
 \ee
where each block represents an $m\times m$ matrix.  The diagonal and super-diagonal blocks represent coefficients that are related to orthogonalization and the subdiagonal blocks represent blocks related to normalization.  The symmetry
tells us that the normalization block $\vek B_{i+1}$, generated when normalizing $\vek V_{i+1}$, can be reused as an orthogonalization coefficient at step $i+1$, where it can be shown that in exact arithmetic
\be\label{eqn.3-term-recur}
	\vek V_{i+1}\vek B_{i+1} = \vek N^{-1}\vek D^{2}\vek V_{i} - \vek V_{i}\vek H_{i}  - \vek V_{i-1}\vek B_{i}.
\ee
We denote $\boldsymbol\CV_{j} = \bbmat \vek V_{1},\vek V_{2},\cdots, \vek V_{j} \ebmat\in\R^{N\times jm}$ to be the matrix whose columns are the block Lanczos vectors.  
In terms of the matrix $\vek T_{j}$, we can use \cref{eqn.3-term-recur} to derive the Lanczos relation, relating $\vek N^{-1}\vek D^{2}$, 
$ \boldsymbol\CV_{j}$, and $\underline{\vek T}_{j}$,
\be\label{eqn.lanczos-relation}
	\vek N^{-1}\vek D^{2} \boldsymbol\CV_{j} =  \boldsymbol\CV_{j+1}\underline{\vek T}_{j}.
\ee
The construction of $\underline{\vek T}_{j}$ is shown as a part of Algorithm \ref{alg.Lanczos}; see, e.g.,  \cite{Saad.Iter.Meth.Sparse.2003,OLeary1980} for further details.
 \RestyleAlgo{boxruled}
\begin{algorithm2e}
\caption{Symmetric Lanczos process for iteratively generating a basis for $\K_{j+1}\prn{\vek N^{-1}\vek D^{2}, \widehat{\vek Y}} $\label{alg.Lanczos}}
Given: starting block of vectors $\widehat{\vek Y} \in\R^{N\times m}$ \\
Out:  $\prn{\cdot,\cdot}_{\vek N}$-orthonormal basis $\curl{\vek V_{1},\vek V_{2},\ldots,\vek V_{j+1}}$\\
Orthogonalize starting block in $\prn{\cdot,\cdot}_{\vek N}$ producing {\tt QR}-factorization: $\widehat{\vek Y}= \underbrace{\vek V_{1}}_{\in\R^{N\times m}}\underbrace{\vek B_{0}}_{\in\R^{m\times m}}$
with $\vek V_{1}^{T}\vek N\vek V_{1} = \vek I_{m}$\\
\For{$i = 1,2,\ldots, j$}{
	$\vek W \leftarrow \vek N^{-1}\vek D^{2}\vek V_{i}$\\
	\If{$i>1$}{
		$\vek W \leftarrow \vek W - \vek V_{i-1}\vek B_{i}$
	}
	$\vek H_{i} = \vek V_{i}^{T}\vek N\vek W$\\
	$\vek W\leftarrow \vek W - \vek V_{i}\vek H_{i}$\\
	Orthogonalize $\vek W$ in $\prn{\cdot,\cdot}_{\vek N}$ producing {\tt QR}-factorization: $\vek W = \vek V_{i+1}\vek B_{i+1}$
}
\end{algorithm2e}
The symmetric Lanczos relation \eqref{eqn.lanczos-relation}
follows directly from the setup of this process.  Furthermore, since the columns of $\boldsymbol\CV_{j}$ 
form an $\vek N$ orthonormal basis for $\K_{j}\prn{\vek N^{-1}\vek D^{2},\widehat{\vek Y}}$, we have that $ \boldsymbol\CV_{j}^{T}\vek N \boldsymbol\CV_{j+1} = \bbmat \vek I_{jm} & \vek 0 \ebmat$, and 
one can then show that,
\begin{align}
	 \vek N^{-1}\vek D^{2} \boldsymbol\CV_{j} &=  \boldsymbol\CV_{j+1}\underline{\vek T}_{j}\nn\\
	 \boldsymbol\CV_{j}^{T}\vek N \vek N^{-1}\vek D^{2} \boldsymbol\CV_{j} &=  \boldsymbol\CV_{j}^{T}\vek N\boldsymbol\CV_{j+1}\underline{\vek T}_{j}\nn\\
	 \boldsymbol\CV_{j}^{T}\vek D^{2} \boldsymbol\CV_{j} &=  \bbmat \vek I_{jm} & \vek 0 \ebmat\underline{\vek T}_{j}\nn\\
	 \boldsymbol\CV_{j}^{T}\vek D^{2} \boldsymbol\CV_{j} &=  \vek T_{j},\label{eqn.tridiag-factorization}
\end{align}
which means that $\vek T_{j}$ is symmetric. 

Building the basis in blocks with the computational kernel being a matrix-times-matrix product and also block orthogonalization produces an
algorithm with computationally attractive properties in terms of computational intensity (i.e., the amount of calculations done per unit of data moved into and out of cache) \cite{BDJ.2006,Parks.Soodhalter.Szyld.16}, and this is a topic of current study in the context of GPU implementation \cite{HigginsGPUTalk}.  Unlike the classical CG approach, though, it is non-trivial to develop a short-recurrence approach to treat symmetric Sylvester equations.  However, an effective, efficient approach has been described in the literature. 

% %%%%%%%%%%%%%%%%%%%%%%%%%%%%%%%%%%%%%%%%%%%%%%%%%%%%%%%%%%%%%%%%%%%%%%%%%%%%%%%%%%%%%%%%%%%%%%%%%%%%%%%%%%%%%%%%%%%%%%%%%%
\subsubsection{Projection methods for Sylvester equations}
Simoncini and Palitta developed a solver for large-scale Sylvester problems where both matrices are symmetric \cite[Section 5.2]{Palitta2018a}.
We adapt this work  to develop a practical method for treating \cref{eqn.Axb-Sylv}.
Similar to methods for solving linear systems, we can discretize the Sylvester equations by taking approximations from a subspace and solving a constrained set of equations to determine the specific approximation.
In this case, the unknown from \cref{eqn.Axb-Sylv} is $\vek M\in\R^{N\times m}$.  At step $j$, we make an approximation of the form 
\be\label{eqn.soln-darstellung}
	\vek M\approx\vek M_{j} = \boldsymbol\CV_{j}\vek Z_{j} = \sum_{i=1}^{j}\vek V_{i}\vek G_{i}  \mwhere \vek G_{i}\in\R^{m\times m}\mand\vek Z_{j} = \bbmat \vek G_{1}\\ \vek G_{2}\\ \vdots \\ \vek G_{j} \ebmat\in\R^{jm\times m}.
\ee
Substituting into \cref{eqn.Axb-Sylv} yields
\be\nn
	\vek N^{-1}\vek D^{2}\boldsymbol\CV_{j}\vek Z_{j} + \boldsymbol\CV_{j}\vek Z_{j}\widehat{\vek L}  \approx \widehat{\vek Y}.
\ee
We determine $\vek Z_{j}$ to be the solution satisfying the equation
\be\label{eqn.Sylv-weak-formulation}
	\boldsymbol\CV_{j}^{T}\vek N\prn{\vek N^{-1}\vek D^{2}\boldsymbol\CV_{j}\vek Z_{j} + \boldsymbol\CV_{j}\vek Z_{j}\widehat{\vek L}}  = \boldsymbol\CV_{j}^{T}\vek N\widehat{\vek Y}.
\ee
Let $\vek E_{1}\in\R^{jm\times m}$ be the matrix with the identity matrix in the first $m\times m$ block and zeros elsewhere.
Observing that $\widehat{\vek Y} = \vek V_{1}\vek B_{0} = \boldsymbol\CV_{j}\vek E_{1}\vek B_{0}$ and using \cref{eqn.tridiag-factorization}, we can simplify \cref{eqn.Sylv-weak-formulation} as
\be\label{eqn.Sylv-small}
	\vek T_{j}\vek Z_{j} + \vek Z_{j}\widehat{\vek L} = \vek E_{1}\vek B_{0}.
\ee	
The task has now become to solve a small, fairly dense set of Sylvester equations.  We can use standard methods to treat this problem. 

We approach \cref{eqn.Sylv-small} by decomposing the matrices $\vek T_{j}$ and $\widehat{\vek L}$.  The matrix $\vek T_{j}$ is symmetric, block-tridiagonal.  Thus, we can efficiently compute the eigendecomposition
$\vek T_{j} = \vek Q_{j}\boldsymbol\Phi_{j}\vek Q_{j}^{T}$ where $\vek Q_{j}^{T}\vek Q_{j}=\vek I_{j}$ and $\boldsymbol\Phi_{j}=\diag\curl{\phi_{1},\phi_{2},\ldots,\phi_{j}}$ is the diagonal matrix of eigenvalues of $\vek T_{j}$.  The matrix $\widehat{\vek L}$ is not symmetric with respect
to the standard Euclidean inner product.  However, we prove it is symmetric with respect to another inner product.
\blem
	The matrix $\widehat{\vek L}=\vek A^{T}\vek T\vek A\vek P^{-1}$ is symmetric with respect to the inner product induced by $\vek P^{-1}$.
\elem
\bproof
	The matrix $\vek P$ is a diagonal matrix with positive entries on the diagonal.  Thus it and its inverse are symmetric positive-definite, and $\vek P^{-1}$ induces the inner product $\prn{\cdot,\cdot}_{\vek P^{-1}}$, with
	$\prn{\vek x,\vek y}_{\vek P^{-1}} = \vek y^{T}\vek P^{-1}\vek x$.  We can show by direct calculation that $\widehat{\vek L}$ is symmetric with respect to $\prn{\cdot,\cdot}_{\vek P^{-1}}$,
	\be\nn
		\prn{\vek A^{T}\vek T\vek A\vek P^{-1}\vek x,\vek y}_{\vek P^{-1}} = \vek y^{T}\vek P^{-1}\vek A^{T}\vek T\vek A\vek P^{-1}\vek x=\prn{\vek x,\vek A^{T}\vek T\vek A\vek P^{-1}\vek y}_{\vek P^{-1}}.
	\ee
\eproof
Since $\widehat{\vek L}$ is symmetric with respect to $\prn{\cdot,\cdot}_{\vek P^{-1}}$, it is diagonalizable and its eigenvectors form an orthonormal basis with respect to $\prn{\cdot,\cdot}_{\vek P^{-1}}$.
Thus, we can write
\begin{align*}
    \widehat{\vek L} = \vek U\vek R\vek U^{-1} = \vek U\vek R\vek U^{T}\vek P^{-1}
\end{align*}
with $\vek R=\diag\curl{\rho_{1},\rho_{2},\ldots,\rho_{m}}$ being a diagonal matrix of eigenvalues, and $\vek U^{T}\vek P^{-1}\vek U=\vek I_{m}$.  This decomposition can be computed once at the beginning
of execution.  Inserting this into \cref{eqn.Sylv-small}, we get
\begin{align}
    \nonumber
    \vek Q_{j}\boldsymbol\Phi_{j}\vek Q_{j}^{T}
    \vek Z_{j} 
    + 
    \vek Z_{j}
    \vek U\vek R\vek U^{-1}
	&= 
	\vek E_{1}\vek B_{0}
	\\
	\nonumber
	\iff
    \boldsymbol\Phi_{j}\vek Q_{j}^{T}
    \vek Z_{j} 
    \vek U
    + 
    \vek Q_{j}^T
    \vek Z_{j}
    \vek U\vek R
	&= 
	\vek Q_{j}^T
	\vek E_{1}\vek B_{0}
	\vek U
	\\
	\label{eqn.diag-sylvester}
	\iff
    \boldsymbol\Phi_{j}
    \widehat{\vek Z}_{j} 
    + 
    \widehat{\vek Z}_{j} 
    \vek R
	&= 
	\vek Q_{j}^T
	\vek E_{1}\vek B_{0}
	\vek U,
\end{align}
where $\widehat{\vek Z}_{j} = \vek Q_{j}^T \vek Z_{j}\vek U\iff\vek Z_j = \vek Q_j\widehat{\vek Z}_{j}\vek U^{-1} = \vek Q_j\widehat{\vek Z}_{j}\vek U^{T}\vek P^{-1}$.
Vectorizing \cref{eqn.diag-sylvester}, we get
\begin{align*}
    \prn{\vek R\otimes\vek I + \vek I\otimes \boldsymbol\Phi_{j}}\mathrm{vec}\prn{ \widehat{\vek Z}_{j} } = \mathrm{vec}\prn{\vek Q_{j}^T\vek E_{1}\vek B_{0}\vek U}.
\end{align*}
Observe that $\prn{\vek R\otimes\vek I + \vek I\otimes \boldsymbol\Phi_{j}}$ is a diagonal matrix with diagonal entries of the form $\phi_{i} + \rho_{k}$.
It follows that just as in 
\cite[Section 5.2]{Palitta2018a}, one can directly compute the solution 
\be\label{eqn.soln-constr}
	\vek Z_{j} 
    = 
    \vek Q_{j}\vek F_{j}\vek U^T\vek P^{-1}
    \mwhere 
    \vek F_{j}
    :=
    \prn{\vek F_{ik}} 
    = 
    \frac{\vek e_{i}^{T} \vek Q_{k}^{T}\vek E_{1}\vek B_{0}\vek U \vek e_{k}}{\phi_{i} + \rho_{k}}
    \in
    \R^{j\times m}.
\ee
\brem
	Although $\vek T_{j}$ contains $\vek T_{j-1}$ as its upper-left $(j-1)m\times (j-1)m$ submatrix, the 
	eigendecomposition changes completly and must be recomputed at every iteration.
\erem

To make this approach computationally efficient, we must forego computing the solution approximations during the iteration.  Rather, we run this iteration and use a cheap residual norm approximation to monitor convergence.  In \cite[Section 5.2]{Palitta2018a}, it is shown 
how to calculate an expresson for the residual norm without forming the solution explicitely. 
This done via accumulation, and in \cite[Algorithm 5]{Palitta2018a}, we see that, in principle, one needs
to calculate the eigendecomposition of $\vek T_{i}$ at each iteration $i$.  However, the authors discuss
in \cite[Section 3.2]{Palitta2018a} how one can calculate only the needed entries of the eigenvectors to 
calculate the residual norm. We note that, in our code, we do not take advantage of the trick and simply compute
the full eigendecomposition.  This is what we show in our description of the procedure, 
Algorithm \ref{alg.sym-sylv-solver}.

Once the residual norm estimate reaches tolerance, the iteration can be stopped, and we construct the approximate solution. We construct $\vek Z_{j}$ and rerun the Lanczos process to get the vectors
again for constructing the solution, one block Lanczos vector at a time.  
% Computing the solution and checking the residual error at each iteration would be costly.  Fortunately, we can 
% compute the residual norm at each iteration without constructing the solution.  Thus, we can run the Lanczos 
% iteration to build $\vek T_{j}$, use this and other quantities to check the residual norm.  Once the convergence
% criteria has been satisfied, we stop, construct $\vek Z_{j}$, and rerun the Lanczos process to get the vectors
% again for constructing the solution. In \cite[Section 5.2]{Palitta2018a}, it is shown 
% how to calculate an expresson for the residual norm without forming the solution explicitely. 
% This done via accumulation, and in \cite[Algorithm 5]{Palitta2018a}, we see that this in principle, one needs
% to calculate the eigendecomposition of $\vek T_{i}$ at each iteration $i$.  However, the authors discuss
% in \cite[Section 3.2]{Palitta2018a} how one can calculate only the needed entries of the eigenvectors to 
% calculate the residual norm.  In this version of the algorithm, we ignore this trick and simply compute
% the full eigendecomposition.  We incorporate this calculation into our description of the procedure, 
% Algorithm \ref{alg.sym-sylv-solver}.
This is because, from \cref{eqn.soln-darstellung}, we can construct the approximation $\vek M_{j}$ block-for-block meaning when we rerun the Lanczos iteration, we can use the blocks as they are constructed and then
	discard them, retaining the memory efficiency of Algorithm \ref{alg.sym-sylv-solver}. 
 \RestyleAlgo{boxruled}
\begin{algorithm2e}
\caption{\label{alg:Sylvester} Symmetric Sparse Sylvester Solver\label{alg.sym-sylv-solver}}
Given: starting block of vectors $\widehat{\vek Y} \in\R^{N\times m}$; coefficient matrix $\vek N^{-1}\vek D^{2}\in\R^{N\times N}$; a small coefficient matrix $\widehat{\vek L}\in\Rmm$; $\varepsilon>0$ a convergence tolerance \\
Out:  An approximation $\vek M\in\R^{N\times m}$ to the solution of the Sylvester equations 
\cref{eqn.Axb-Sylv}\\
Compute eigendecomposition $\widehat{\vek L} = \vek U\vek R\vek U^{T}$\\
Orthogonalize starting block in $\prn{\cdot,\cdot}_{\vek N}$ producing {\tt QR}-factorization: $\widehat{\vek Y}= \vek V\vek B_{0}$\\
$\vek V_{old} = \vek 0$; $i=0$\\
\While{$\sqrt{\gamma} > \varepsilon\nm{\vek B_{0}}_{F}$}{
$i\leftarrow i+1$\\
\tcc{Lanczos step ******************************************}
$\vek W \leftarrow \vek N^{-1}\vek D^{2}\vek V$\\
\If{$i>1$}{
	$\vek W \leftarrow \vek W - \vek V_{old}\vek B_{i}$
}
$\vek H_{i} = \vek V^{T}\vek N\vek W$\\
$\vek W\leftarrow \vek W - \vek V\vek H_{i}$\\
$\vek V_{old}\leftarrow\vek V$\\
Orthogonalize $\vek W$ in $\prn{\cdot,\cdot}_{\vek N}$ producing {\tt QR}-factorization: $\vek W = \vek V\vek B_{i+1}$\\
\tcc{computing residual norm according to \cite[Section 5.2]{Palitta2018a} *********************}
$\gamma\leftarrow 0$\\
Compute eigendecomposition $\vek T_{i}=\vek Q_{i}\boldsymbol\Phi_{i}\vek Q_{i}^{T}$\\
 $\vek S\leftarrow\prn{\vek U^{T}\vek B_{0}^{T}}\prn{\vek E_{1}^{T}\vek Q_{i}}$\\
$\vek J \leftarrow \prn{\vek Q_{i}\vek E_{m}}\vek B_{i+1}^{T}$\\
\For{$j=1,2,\ldots,m$}{
	$\vek K\leftarrow \rho_{j}\vek I + \boldsymbol\Phi_{i}$\\
	$\gamma \leftarrow \gamma + \norm{\prn{\vek e_{j}^{T}\vek S}\vek K^{-1}\vek J}^{2}$
}
}
\tcc{constructing approximation coefficients entry-wise *******}
$\vek F=\vek 0\in\R^{i\times m}$\\
\For{$k=1,2,\ldots, i$}{
	\For{$\ell=1,2,\ldots, m$}{
		$\vek F_{k\ell} = \frac{\vek e_{k}^{T} \vek Q_{\ell}^{T}\vek E_{1}\vek B_{0}\vek U \vek e_{\ell}}{\phi_{k} + \rho_{\ell}}$; see \eqref{eqn.soln-constr}
	}
}
$\vek Z = \vek Q_{i}\vek F\vek U^{T}$\\
\tcc{rerun Lanczos to assemble solution ***********************}
$\vek M = \vek 0$\\
$\vek V \leftarrow \widehat{\vek Y}\vek B_{0}^{-1}$\\
$\vek V_{old} = \vek 0$\\
\For{$j=1,2,\ldots, i$}{
$\vek M \leftarrow \vek M + \vek V\vek G_{j}$, according to \eqref{eqn.soln-darstellung}\\
\tcc{Lanczos step using already-computed coefficients ******}
$\vek W \leftarrow \vek N^{-1}\vek D^{2}\vek V$\\
\If{$j>1$}{
	$\vek W \leftarrow \vek W - \vek V_{old}\vek B_{j}$
}
$\vek W\leftarrow \vek W - \vek V\vek H_{j}$\\
$\vek V_{old}\leftarrow\vek V$\\
$\vek V \leftarrow \vek W\vek B_{j+1}^{-1}$
}
Output $\vek M$
\end{algorithm2e}

\begin{remark}
    In Algorithm \ref{alg:Sylvester} we do not perform any low-rank compression of the block Lanczos vectors.  In our experiments, we determined that the block Lanczos vectors and the residuals generated when running Algorithm \ref{alg:Sylvester} remained numerically full rank (i.e., rank 4).  Thus, we do not pursue a strategy of compression in our Sylvester solver.  
\end{remark}

% %%%%%%%%%%%%%%%%%%%%%%%%%%%%%%%%%%%%%%%%%%%%%%%%%%%%%%%%%%%%%%%%%%%%%%%%%%%%%%%%%%%%%%%%%%%%%%%%%%%%%%%%%%%%%
%\subsection{Solving in implicit Kronecker versus Sylvester form}

% %%%%%%%%%%%%%%%%%%%%%%%%%%%%%%%%%%%%%%%%%%%%%%%%%%%%%%%%%%%%%%%%%%%%%%%%%%%%%%%%%%%%%%%%%%%%%%%%%%%%%%%%%%%%%
% %%%%%%%%%%%%%%%%%%%%%%%%%%%%%%%%%%%%%%%%%%%%%%%%%%%%%%%%%%%%%%%%%%%%%%%%%%%%%%%%%%%%%%%%%%%%%%%%%%%%%%%%%%%%%
\section{Evaluation and Application}
\label{sec:evaluation}
% %%%%%%%%%%%%%%%%%%%%%%%%%%%%%%%%%%%%%%%%%%%%%%%%%%%%%%%%%%%%%%%%%%%%%%%%%%%%%%%%%%%%%%%%%%%%%%%%%%%%%%%%%%%%%
% %%%%%%%%%%%%%%%%%%%%%%%%%%%%%%%%%%%%%%%%%%%%%%%%%%%%%%%%%%%%%%%%%%%%%%%%%%%%%%%%%%%%%%%%%%%%%%%%%%%%%%%%%%%%%

% {\color{red} SPW: Below is what Jiajin did. I think what we need to do here is:
%\begin{itemize}
%    \item ACCURACY: Simulated data where we have the ground truth and can evaluate mean square error. Do these as a function of problem size --- Dung is looking at this now.
 %   \item SCALABILITY: time/storage requirements as a function of problem size, compare to 'straw man' approach of calculating Cholesky --- Jiajin did this below and it's a question of whether we want to repeat it; KMS: also, what is the proper way to test scaling on a GPU? Were this an MPI-based computation, I would suggest investigating the strong- and weak- scaling of the code, but I don't think that makes sense on a GPU
 %   \item PROFILING: profile code to see what proportion of computation time is taken up with different parts of the algorithms, again as a function of problem size --- Dung will look at this next;
 %   \item REAL DATA: illustration with Planck data --- to do.
% \end{itemize} }

We have proposed using two approaches described in the literature for solving this source separation 
problem by exploiting the Kronecker structure of the linear system.  In the first method, we employ conjugate gradients applied
directly to the Kronecker form of the system, and the structure
is exploited to obtain a matrix-free iteration.  In the latter 
method, we similarly exploit the Sylvester formulation of the problem
using a Lanczos-based iterative method, and we are again able to
do this in a matrix-free way.

In this section, we demonstrate with our numerical results that both approaches
are viable for the underlying application problem.
In theory, though, we should expect better performance from the
method based on the Sylvester formulation.  This has been shown in \cite{pal-sim-kroncond}

Three implementations (two variants of the conjugate gradient approach and the Sylvester method) are applied to \textbf{both} \emph{simulated data} of different sizes (allowing for evaluation of the actual error) as well as \emph{actual Planck full sky maps}.  Results are compared in terms of time cost, memory usage, and accuracy. Code profiling explores how these requirements for different aspects of the algorithm change with problem size.

The experiment is designed as follows. The 3 solutions were run with simulated data having different values of $N_{side}$. Since the source $\vek S$ is specified in each simulation, this allows us to compare the results of the algorithms with the true solution of the simulated problem. 
The problem size is specified as the HEALPix level from 1 to 10. Recalling that the solution dimension at level $h$ is $m \times N_{side} =  m \times 12 \times 4^h$, hence increasing the level by 1 increases the dimension of the solution by 4.

For the simulated and Planck data, $m=4$ sources and $n=9$ images were assumed, as discussed in Appendix \ref{app:A}.  Data were simulated by specifying values for the source vector $\boldsymbol{\mathcal{S}}$ (for simplicity, these are generated as independent samples from a normal distribution), and the matrix $\vek A$ and vector $\vek T$ were defined as in Appendix \ref{app:A}) . That allows the data vector $\boldsymbol{\mathcal{Y}}$ to be randomly generated as in Equation \ref{eq:model}.

Code was written in Python and implemented on a cluster with specification described in Table \ref{tab:hpc_spec}.
\begin{table}[h]
\centering
    \begin{tabular}{@{}lll@{}}
        \toprule
        Property & Value \\
        \cmidrule(r){1-1}\cmidrule(lr){2-2}
        Architecture &  64 \\
        Number of Nodes &  100 \\
        Ram per node & 24GB \\
        RAM & 2.4TB \\
        Clock Speed & 2.66GHz \\
        Total number of cores & 1200 \\
        Theoretical Peak Performance & 12.76TF \\
        Operating System & Scientific Linux 6.x \\
        Allocated Memory & 50GB \\
        Allocated Nodes & 2 \\
        Allocated Cores per Node & 12 \\ 
        \bottomrule
        \end{tabular}
        \caption{\label{tab:hpc_spec} Specification of the machine used in the experiments.}
    \end{table}

\subsection{Accuracy}
\begin{figure}[ht!]
    \centering
    \includegraphics[width=0.45\textwidth]{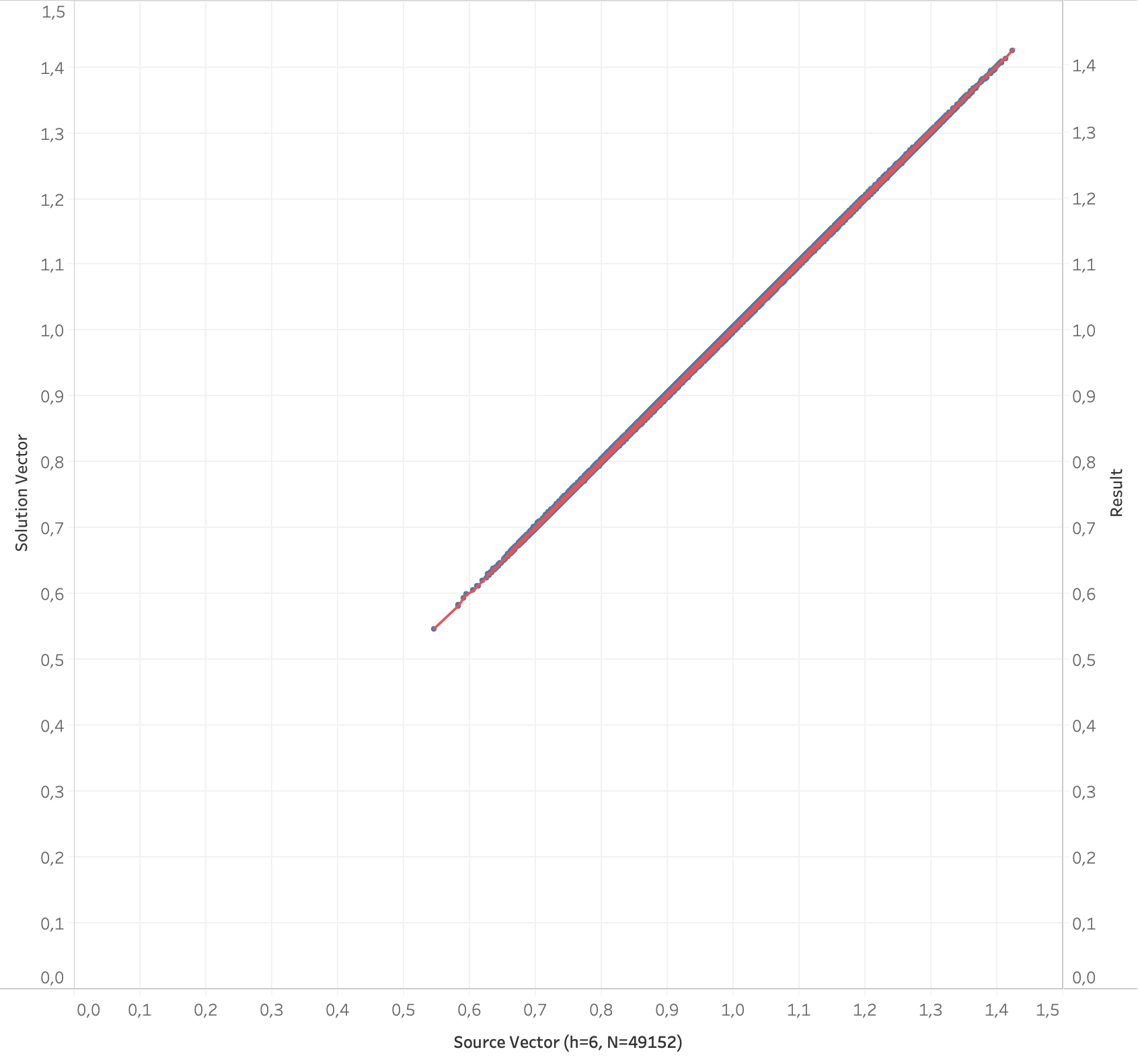}
    \includegraphics[width=0.45\textwidth]{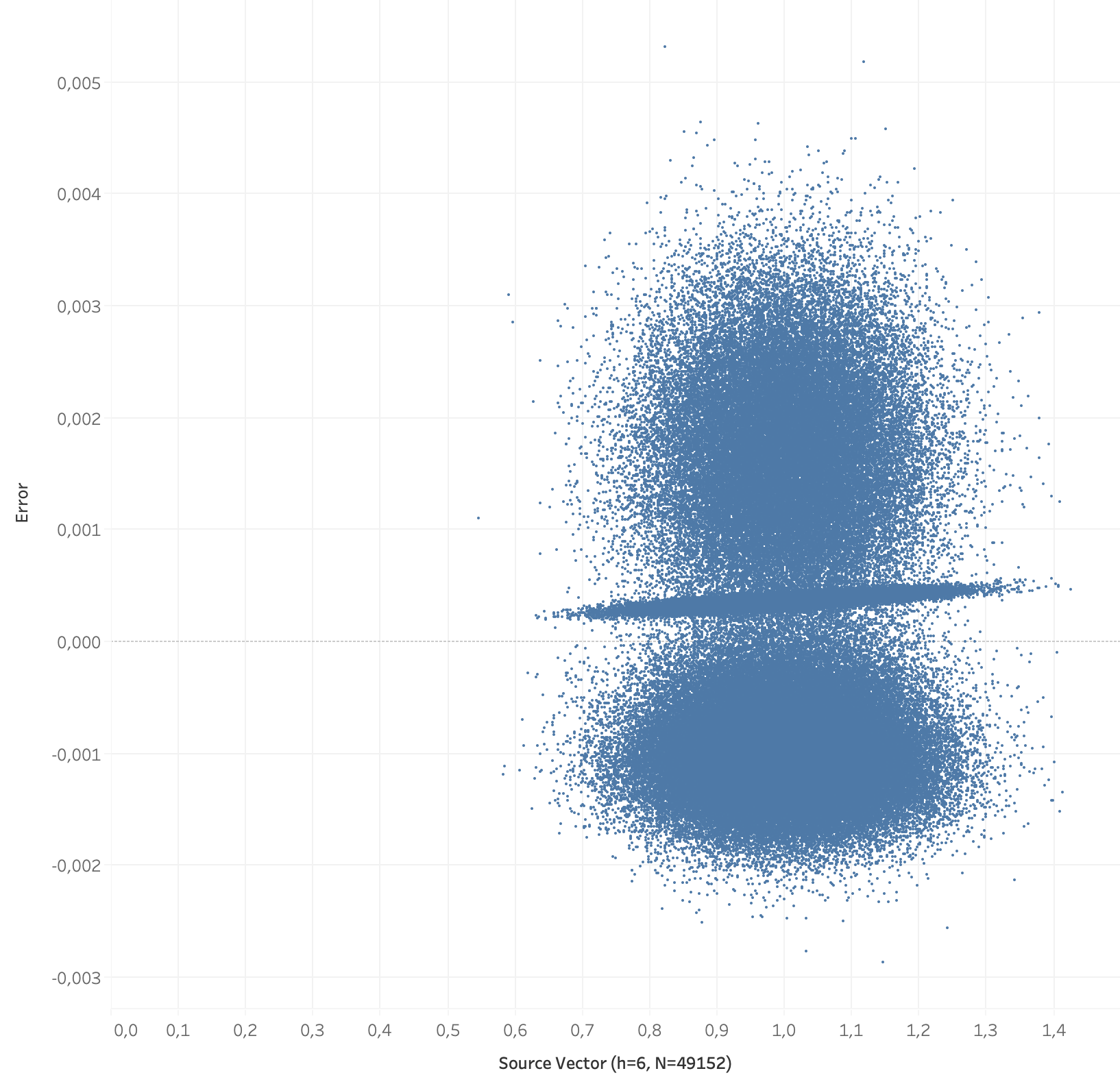}
    \caption{\label{fig:accuracy} Solution accuracy. Scatter plot of components of the true source vector $\mathcal{S}$ versus the solution vector (left) and source vector versus the error (right). }
\end{figure}

Figure \ref{fig:accuracy} shows the relationship between the source vector of the random-generated data and the solution produced by the Conjugate Gradient method on our linear system for simulated data with $h=6$ or $N_{side} = 49,152$. They almost matched  with an error margin (defined to be largest absolute error across all components of the source vector divided by the average value of the source vector) of 0.4\%. Similar results are seen for data of other sizes.

\subsection{Time  and memory usage}
\begin{figure}[!hbp]
  \centering
    \includegraphics[width=0.45\textwidth]{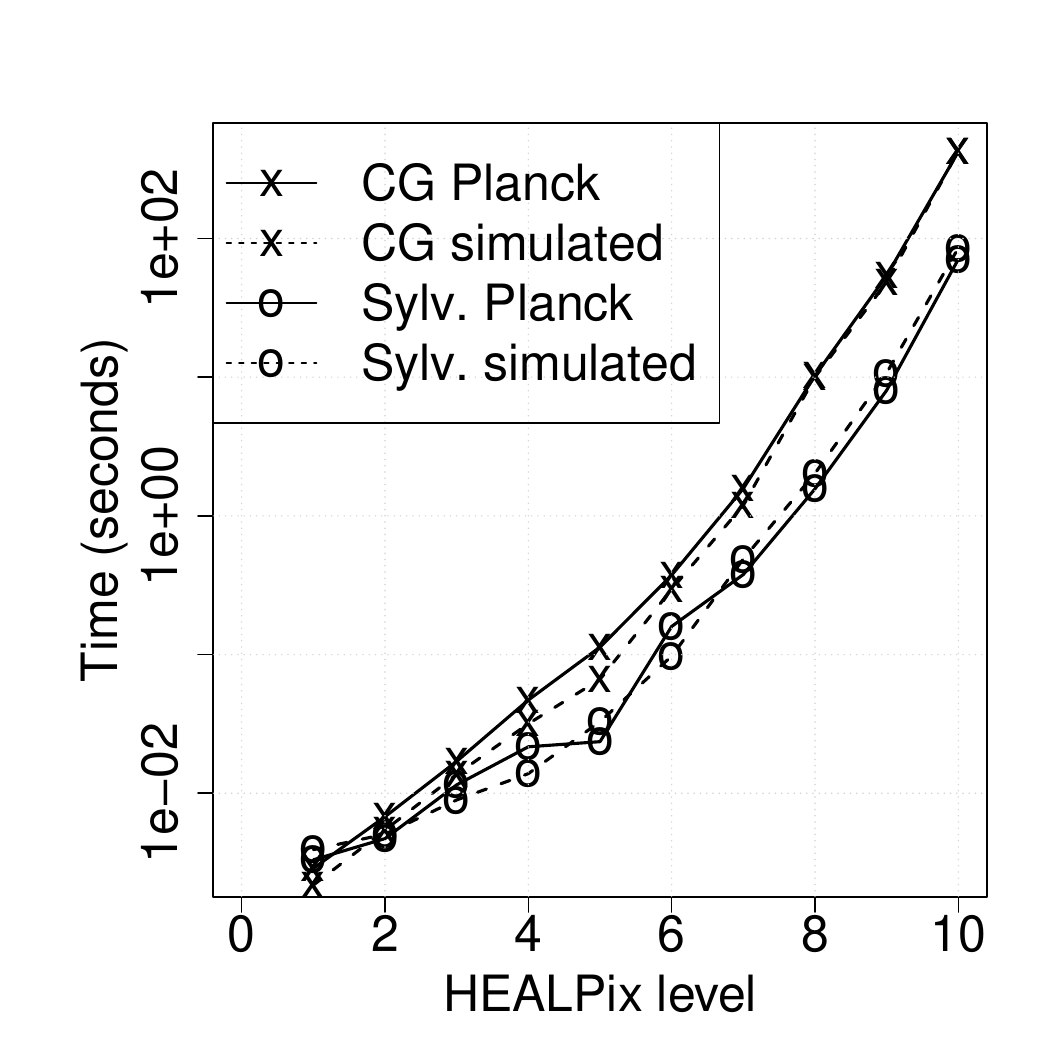}
        \includegraphics[width=0.45\textwidth]{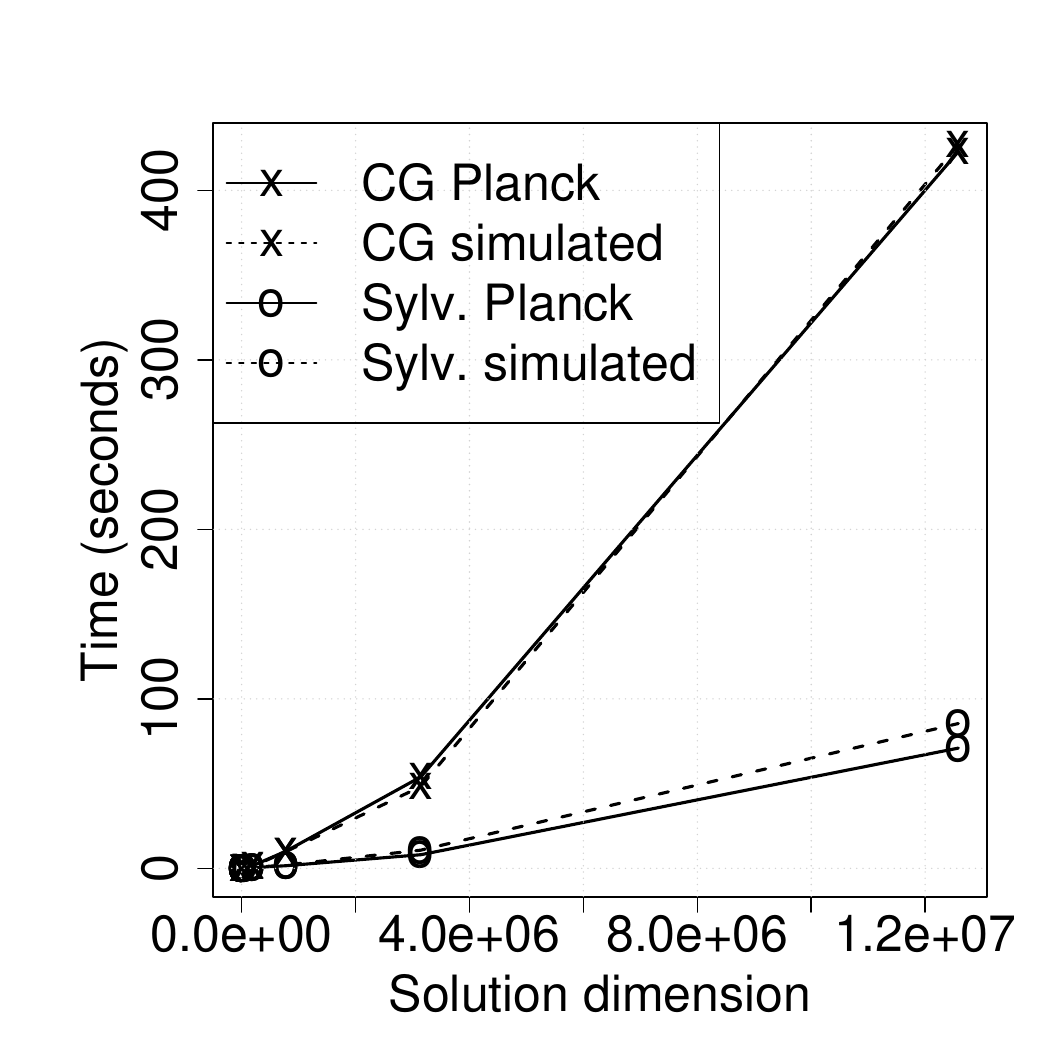}
    \caption{\label{fig:time} Running time (in seconds) by HEALPix level (left) and data size (right).}
\end{figure}

% Table \ref{tab:time} and 
Figure \ref{fig:time} shows the time cost (in seconds) while running the Conjugate Gradient and Sylvester methods on simulated and Planck data, as a function of HEALPix level and the data size; Appendix \ref{app:healpix} explains that HEALPix level $h$ means a data size of $4^h$, hence signifies log data size. The difference between Planck and simulated data is due to the longer time that it took to load in the Planck data from a file than it took to simulate data of the same size. However, the main difference is that the \emph{Sylvester method is consistently faster, taking only about 15\% of the time of the Conjugate Gradient approach for the larger data sizes}.  In both cases, the algorithm demonstrates a close to linear complexity with the data size $N$, a substantial improvement over exact solution methods that are at least quadratic with $N$. We hypothesize that the superior performance of the Sylvester block Krylov method over the CG method applied directly to the Kronecker-structured problem is due both to conditioning issues reported in \cite{pal-sim-kroncond} as well as the superior data movement characteristics (higher computational intensity) often exhibited by block-Krylov-subspace-based methods; see, e.g., \cite{BDJ.2006}.

\begin{figure}[!hbp]
  \centering
    \includegraphics[width=0.45\textwidth]{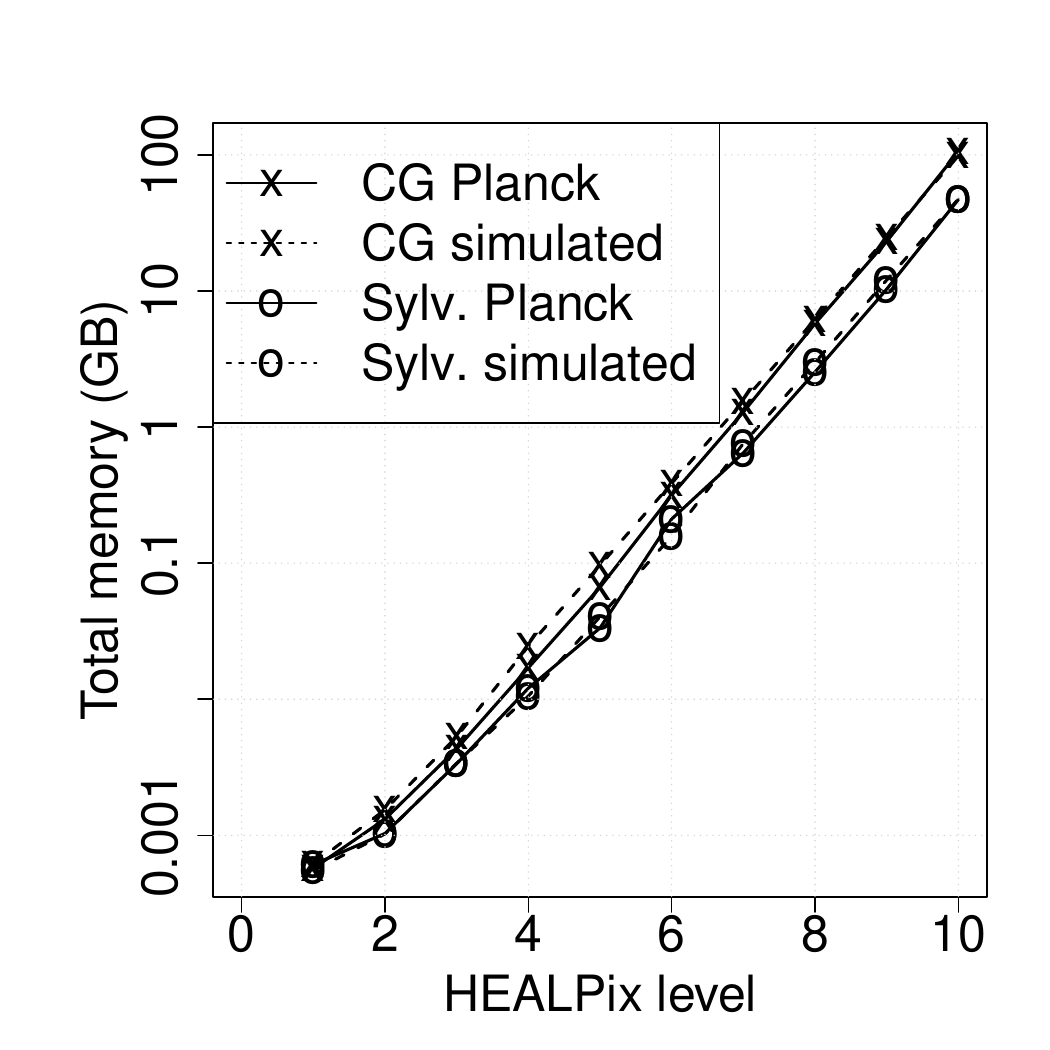}
    \includegraphics[width=0.45\textwidth]{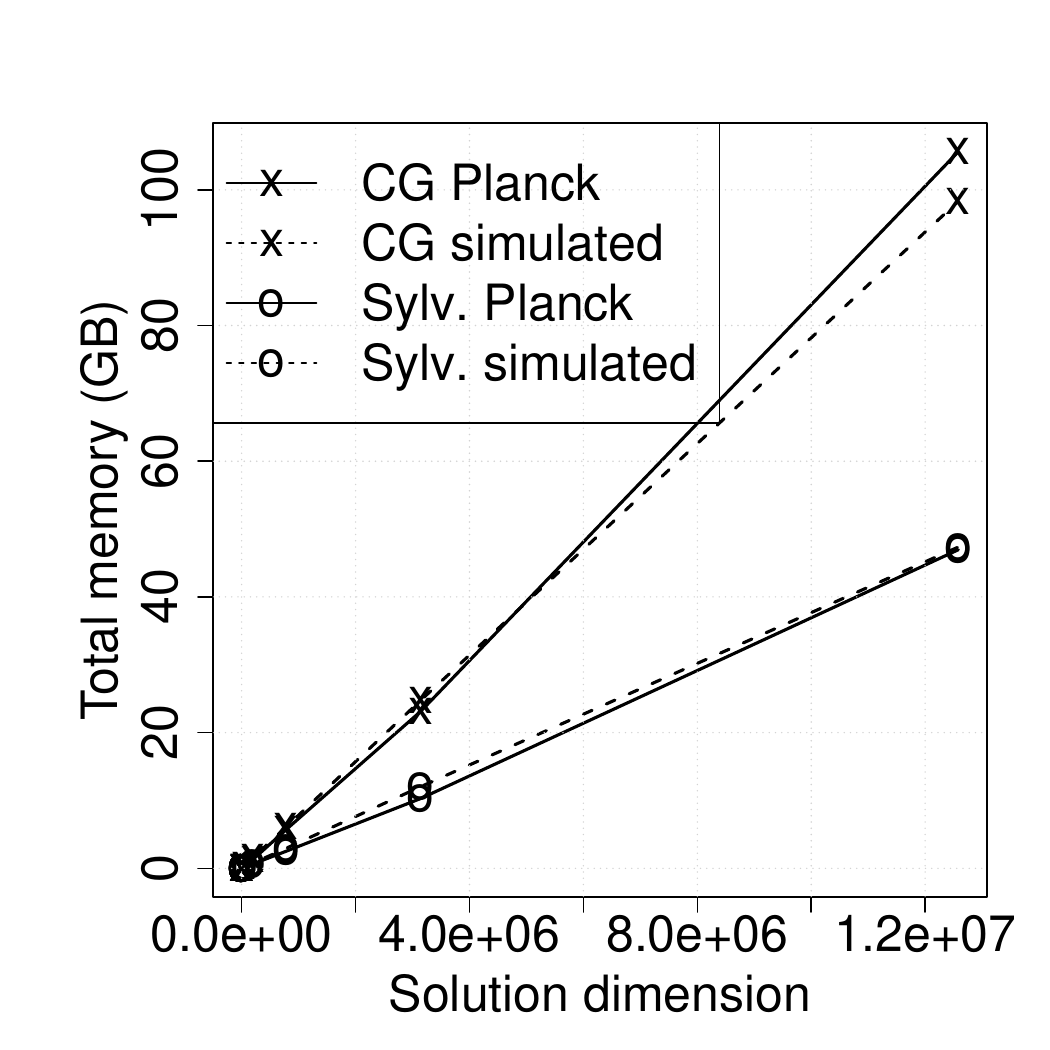}
        \caption{\label{fig:memory}A comparison of total memory usage by HEALPix level (left) and data size (right) for Algorithm \ref{alg.CGnoprecond} (CG) versus Algorithm \ref{alg:Sylvester} (Lanczos-Sylvester -- i.e., Sylv.) for both synthetic data and actual Planck satellite data.}
\end{figure}

Figure \ref{fig:memory} show the total memory usage (in MB) while running the Conjugate Gradient and Sylvester methods on simulated and Planck data. We see that the Sylvester method's memory requirements are about half of that of the Conjugate Gradient, the former needing about 47GB for HEALPix level 10, which is the size of the Planck satellite data. 

% \begin{table}[h]
%\centering
% \begin{tabular}{@{}lll@{}}
%         \toprule
%         HEALPix Level & Random Data & Planck Data \\
%         \cmidrule(r){1-1}\cmidrule(lr){2-2}\cmidrule(l){3-3}
%         1  & 0.094 & 0.131 \\
%         2  & 0.191 & 0.173 \\
%         3  & 0.444 & 0.469 \\
%         4  & 1.535 & 1.662 \\
%         5  & 5.932 & 6.453 \\
%         6  & 26.46 & 25.77 \\
%         7  & 110.8 & 111.7 \\
%         8  & 443.6 & 456.9 \\
%         9  & 1820  & 2025 \\
%         10 & 7281  & 6940 \\
%         \bottomrule
%         \end{tabular}
%     \caption{\label{tab:time} Time cost (in seconds) of the Conjugate Gradient solutions.}
% \end{table}

% \begin{table}[h]
% \centering
% \begin{tabular}{@{}lll@{}}
%         \toprule
%         HEALPix Level & Random Data & Planck Data \\
%         \cmidrule(r){1-1}\cmidrule(lr){2-2}\cmidrule(l){3-3}
%         1 & 0.53 & 1.34 \\
%         2 & 0.46 & 1.44 \\
%         3 & 1.02 & 2.63 \\
%         4 & 3.29 & 7.35 \\
%         5 & 12.35 & 26.43 \\
%         6 & 48.6 & 102.31 \\
%         7 & 193.6 & 408.56 \\
%         8 & 773.59 & 1625.62 \\
%         9 & 3093.54 & 6539.69 \\
%         10 & 12373.58 & 16730.00 \\
%         \bottomrule
%         \end{tabular}
%     \caption{\label{tab:memory} Memory Usage (in MB) of the Conjugate Gradient Solutions}
% \end{table}

\subsection{Code profiling}
\begin{figure}
    \centering
     \includegraphics[width=0.45\textwidth]{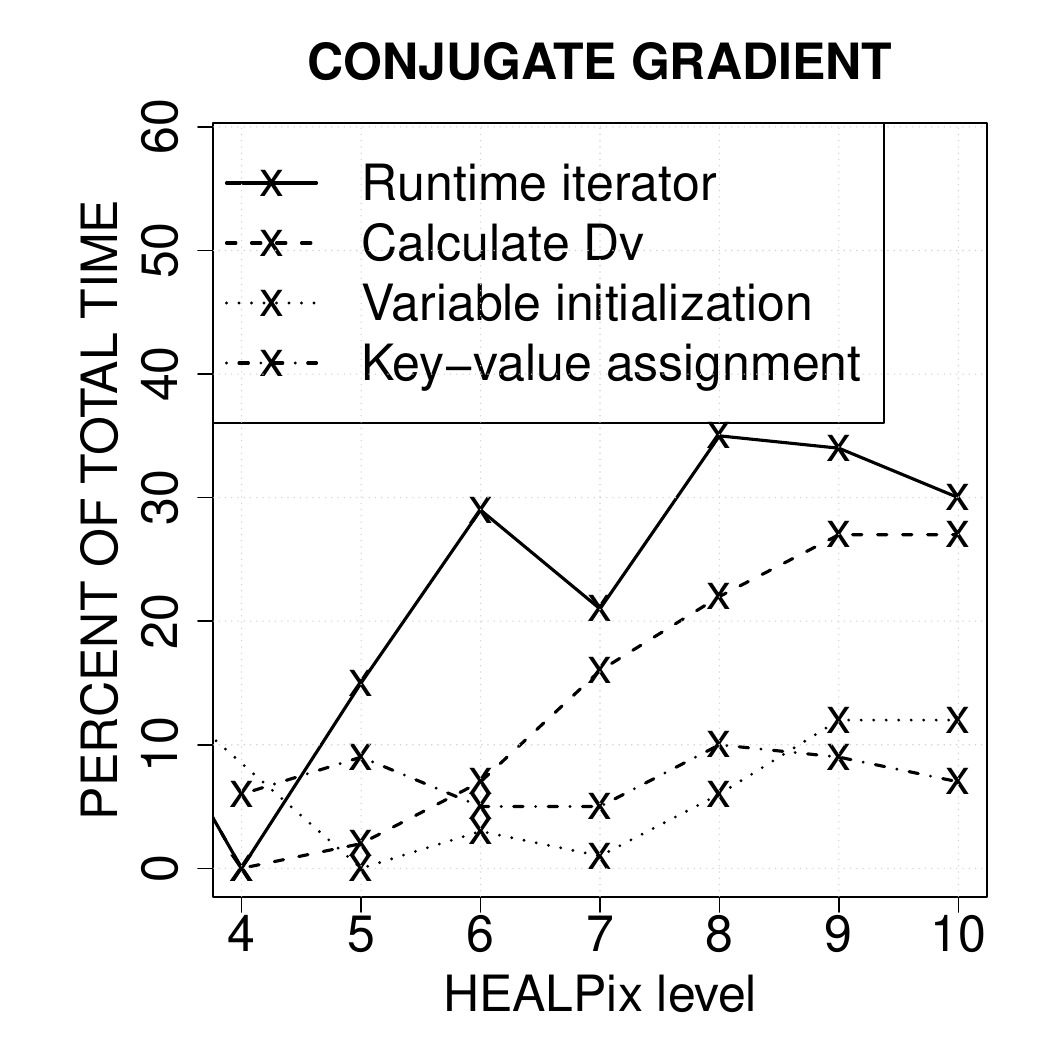} 
     \includegraphics[width=0.45\textwidth]{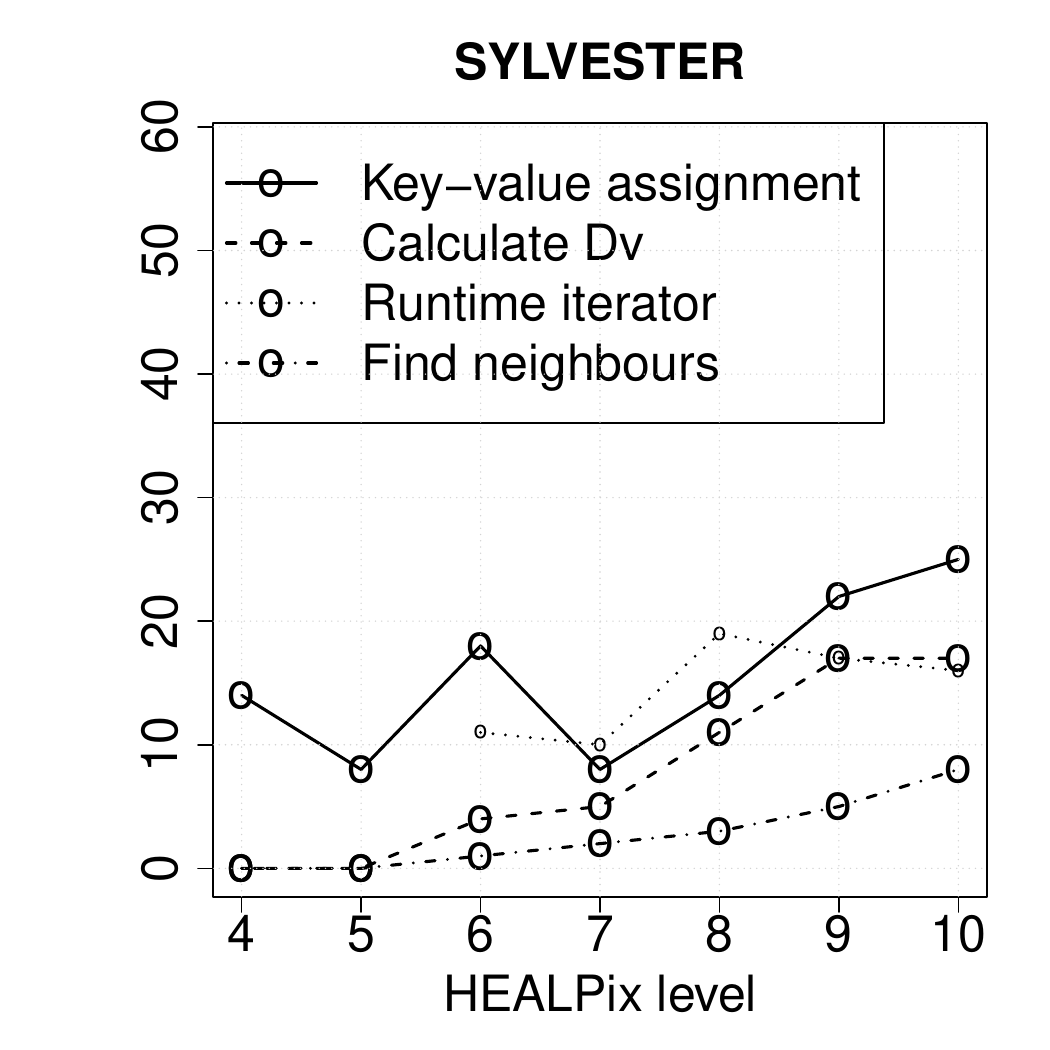} \\  
     \includegraphics[width=0.45\textwidth]{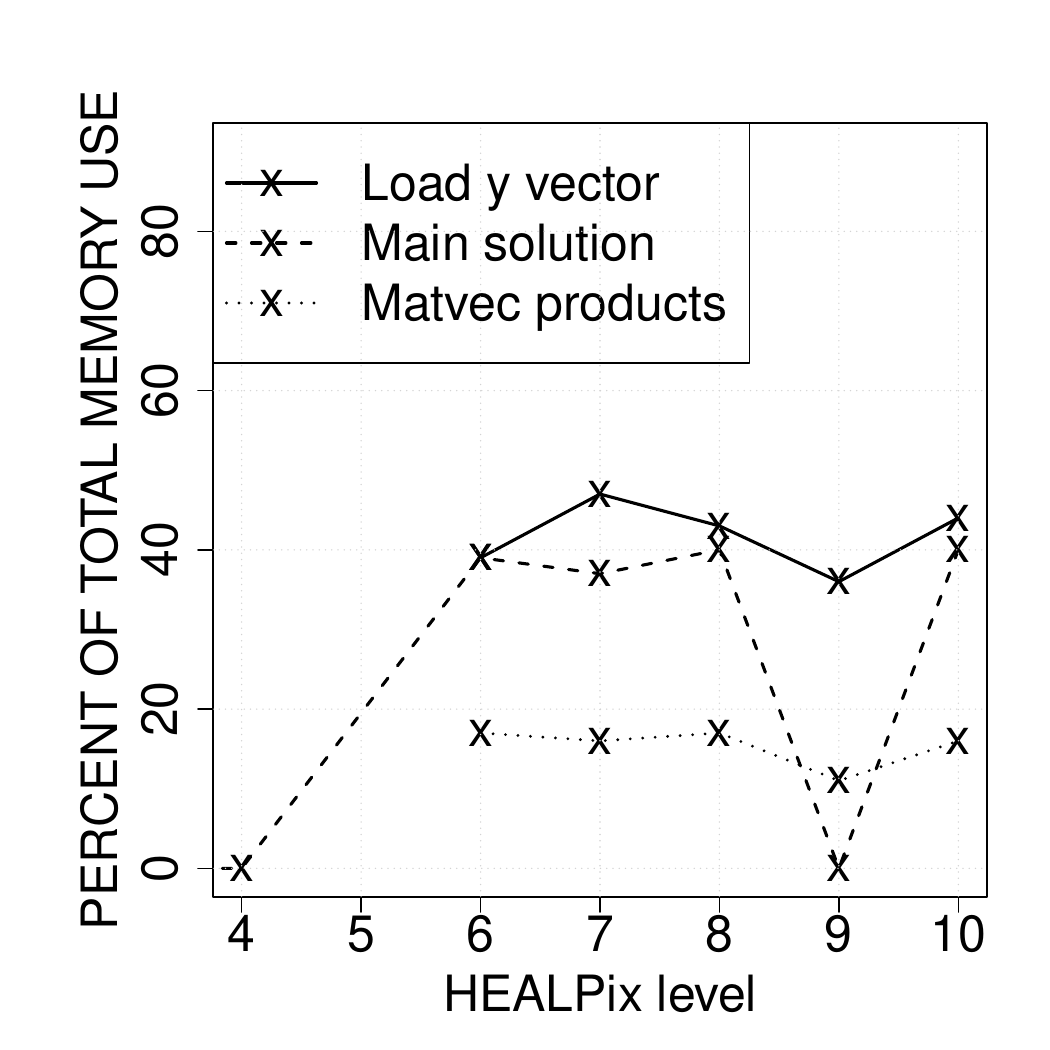}  
     \includegraphics[width=0.45\textwidth]{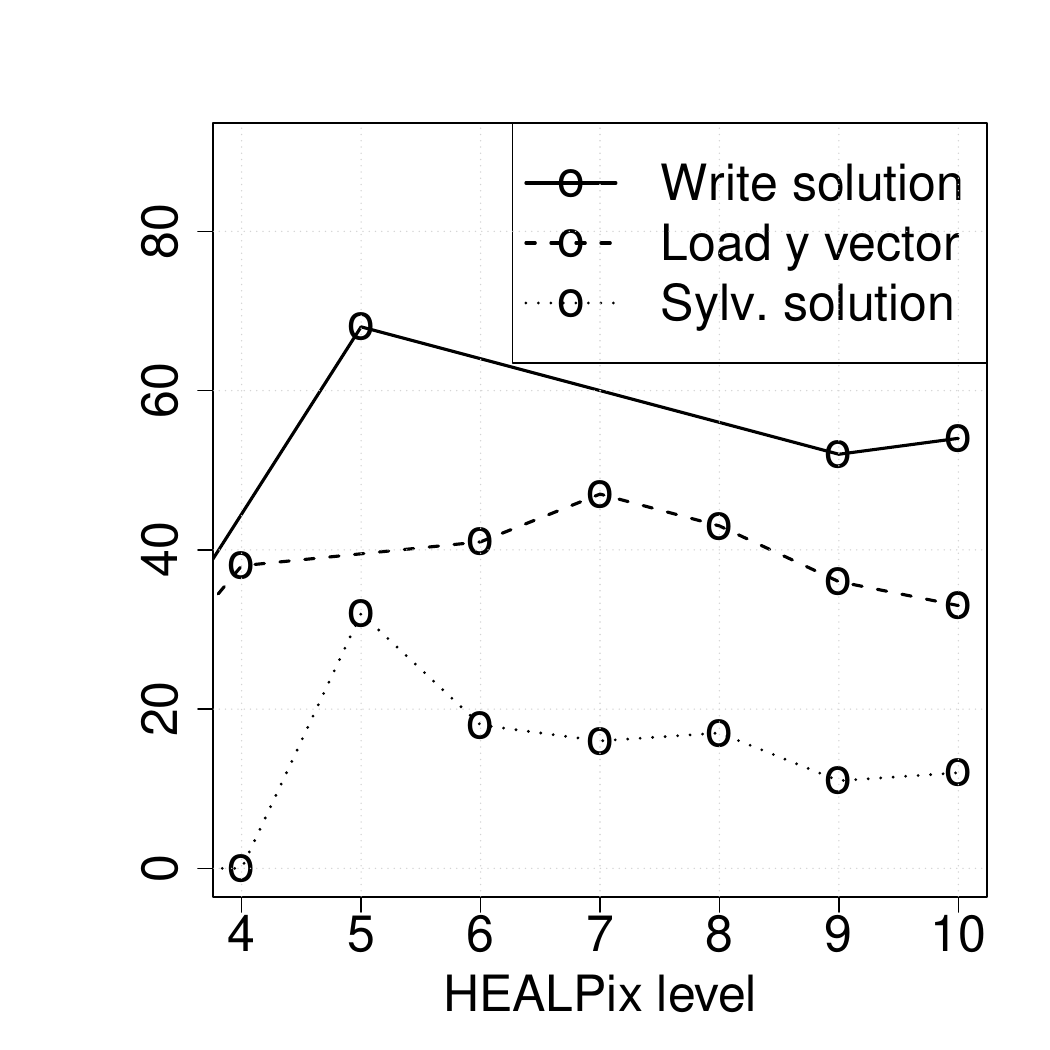} 
    \caption{Code profiling. The  functions that consume the most resources  as a function of problem size. The top row is percentage of total computation time and the bottom row is percentage of maximum memory use. Table \ref{tab:profile_defs} provides more information on the functions. \emph{Sylv.} refers to Algorithm \ref{alg:Sylvester}.}
    \label{fig:profile}
\end{figure}

Figure \ref{fig:profile} is a break down of the computation time and memory use, as a function of problem size, by the functions in the code that used the most of these resources.  A description of these functions is given in \Cref{tab:profile_defs}.   For time, the top 4 functions at HEALPix level 10 are plotted, while for memory there were only 3 functions that accounted for almost all of the memory use. For memory, unsurprisingly it is the functions that load or make use of all of the data that dominate resources. For time, we see that mapping and computation around the $\vek D$ matrix dominate the computation time. This suggests that efforts to further speed up the algorithm should focus on this task.  Encouragingly, operations around the $\vek D$ matrix are well suited to parallelisation \cite{saad89, burr17}.

% The time taken to compute the solution was dominated by matrix-vector multiplications involving the $Q$ matrix, as in Section \ref{subsec:matvecQ}; typically 98\% of the total solution time was spent on this operation. 

\begin{table}
\centering
   \begin{tabular}{@{}lll@{}}
        \toprule
        Function & Description \\
        \cmidrule(r){1-1}\cmidrule(lr){2-2}
        Calculate Dv  & Implementation of Algorithm 1.\\
        Find neigbours &  Identify neighbouring pixel indices in the HEALPix partition. \\
        Runtime iterator & A system function to determine the next item to process, used to parallelise the code. \\
        Key-value assignment & A system function to assign a value to a key (string) in the code. \\
        Variable initialization & Define all variables for the problem and allocate memory. \\
        \cmidrule(r){1-1}\cmidrule(lr){2-2}
        Load $y$ vector & Self explanatory. \\
        Main solution & Managing the main loop of code.\\
        Matvec products & Implementation of the calculations described in Section \ref{subsec:matvec}.\\
        Sylvester solution & Main while loop  of Algorithm \ref{alg:Sylvester}. \\
        Write solution & Self explanatory. \\
        \bottomrule
        \end{tabular}
\caption{\label{tab:profile_defs} Definition of functions listed in Figure \ref{fig:profile}.}
\end{table}
\subsection{Comparison with approach from \cite{kohller}}\label{section:comparison-Kohller}

We construct the sparse-dense form of the Sylvester equation and apply Algorithm \ref{alg:big-small-Kohler} to the Planck dataset. For this dataset, $m = 4$ and $N = 12 \times 4^{h}, h = 1,2,...,10$. For the largest size of the problem, $\mathbf{H} \in \mathbb{R} ^ {12,582,912 \times 12,582,912}$, wherein we remind the reader that this matrix is the large, sparse matrix appearing in \eqref{eq:sparse-dense-Sylvester}.

The results in terms of accuracy and running time of our Sylvester approach Algorithm \ref{alg:Sylvester} and the approach from \cite{kohller}, (i.e., Algorithm \ref{alg:big-small-Kohler}) are shown in \Cref{table:sparse-dense-comparisons} a and b, respectively. The Sparse-Dense method shows significant advantages in accuracy with a comparable running time. However, we observed that for larger Planck data sizes, the solve time for Algorithm \ref{alg:big-small-Kohler} grew quadratically.  Further investigation revealed that unlike for the problems discussed in \cite{kohller}, the SuperLU inner linear solver was not effective at reducing solve time growth rate.    
We note that we did try different sparse direct solvers (UMFPACK, SuperLU) and an iterative method (Conjugate Gradient), all within the scipy Python library. The SuperLU solver proved to be the most \emph{efficient} and accurate one for our specific data.
We tried multiple approaches with this code to reduce the growth in run-time, but the growth in runtime remained quadratic. Further investigation of this issue is out of the scope of our paper. In comparison, our method is demonstrated to be efficient, even for the full size of the Planck dataset.  

\captionsetup[sub]{font=normal,labelfont={bf,sf}}
\begin{table}
\begin{center}
\begin{tabular}{||c c c c||} 
 \hline
Level ($h$) & N & Lanczos-Sylvester & Sparse-Dense Method  \\ [0.5ex] 
 \hline\hline
 1 & $48$ & $1.72\mathrm{e}{-06}$ & $2.23\mathrm{e}{-10}$\\ 
 \hline
 2 & $192$ & $4.21\mathrm{e}{-06}$ & $3.28\mathrm{e}{-10}$\\
 \hline
 3 & $768$ & $7.02\mathrm{e}{-06}$ & $4.73\mathrm{e}{-10}$\\
 \hline
 4 & $3072$ & $4.86\mathrm{e}{-06}$ & $6.19\mathrm{e}{-10}$ \\
 \hline
 5 & $12288$ & $8.62\mathrm{e}{-06}$ & $5.85\mathrm{e}{-09}$\\
 \hline
 6 & $49152$ & $1.08\mathrm{e}{-06}$ & $2.19\mathrm{e}{-10}$\\
 \hline
 7 & $196608$ & $1.35\mathrm{e}{-06}$  & $2.6\mathrm{e}{-09}$ \\
 \hline
 8 & $786432$ & $1.58\mathrm{e}{-06}$ & $1.93\mathrm{e}{-09}$\\
 \hline
 9 & $3145728$ & $1.74\mathrm{e}{-06}$ & $1.21\mathrm{e}{-09}$\\
 \hline
 10 & $12582912$ & $4.23\mathrm{e}{-06}$& $-$ \\
 \hline
\end{tabular} \\
\subcaption{ Residual Comparison}
\end{center}
% \end{table}

% \begin{table}
\begin{center}
\begin{tabular}{||c c c c||} 
 \hline
Level ($h$) & N & Lanczos-Sylvester & Sparse-Dense Method  \\ [0.5ex] 
 \hline\hline
 1 & $48$ & $0.003$ & $0.001$\\ 
 \hline
 2 & $192$ & $0.004$ & $0.006$\\
 \hline
 3 & $768$ & $0.012$ & $0.035$\\
 \hline
 4 & $3072$ & $0.021$ & $0.3$ \\
 \hline
 5 & $12288$ & $0.023$ & $0.755$\\
 \hline
 6 & $49152$ & $0.159$ & $4.35$\\
 \hline
 7 & $196608$ & $0.377$ & $37.312$ \\
 \hline
 8 & $786432$ & $1.566$ & $301.698$\\
 \hline
 9 & $3145728$ & $8.042$ & $6948.4228$\\
 \hline
 10 & $12582912$ & $70.930$ & $-$ \\
 \hline
\end{tabular} \\
\subcaption{ Running Time Comparison}%, Planck Data\label{table:runtime-comp}}
\end{center}
\caption{Comparison using Planck Data (with $m=4$ sources), between the sparse-dense method proposed in \cite{kohller} and the Lanczos-Sylvester method adapted from \cite[Section 5.2]{Palitta2018a}.  The Sparse-Dense method shows significant advantages in accuracy with a comparable running time for smaller problems. However, for larger problems, the solve-time required by the sparse-dense method, even when implemented using Super-LU as the inner solver grows like $\CO(N^2)$ rather than $\CO(N)$ for these problems.
This causes the sparse-dense solver to run so long that we aborted the $h=10$ experiment.
\label{table:sparse-dense-comparisons}}
\end{table}

% \begin{figure}[!hbp]
%   \centering
%     \includegraphics[width=0.45\textwidth]{profiling_Planck_h.pdf}
%     \includegraphics[width=0.45\textwidth]{profiling_sim_h.pdf}
%     \caption{\label{fig:profile_old}Computation time profile by HEALPix level.  On the left, data input and solution for Planck data. On the right, data generation and solution for simulated data.}
% \end{figure}

% Each time, these three methods receive the same input parameters:$m=9$,$n=4$, and $\vek A$ amd $\vek T$ as defined in Appendix \ref{app:A},
% = 
% \begin{bmatrix}
% 1 & 24.31408176 & 0.1808804 & 13.15799486 \\
% 1 & 8.81745815 & 0.3154685 & 5.80101879 \\
% 1 & 2.58070139 & 0.6121030 & 2.15148219 \\
% 1 & 1.00587794 & 1.0058779 & 1.00587794 \\
% 1 & 0.39224839 & 1.6297015 & 0.47074157 \\
% 1 & 0.13192930 & 2.7832105 & 0.19585548 \\
% 1 & 0.03801049 & 4.9311851 & 0.07232168 \\
% 1 & 0.01327148 & 7.7040185 & 0.03151240 \\
% 1 & 0.00509560 & 11.3366157 & 0.01524107 
% \end{bmatrix}
% $ \\
% and   
% \begin{equation}\nonumber
%     \vek T = diag\{629881.6, 694444.4, 783146.7, 12755102.0, 30864197.5, 30864197.5, 30864197.5, 30864197.5, 30864197.5\}.
% \end{equation}

\section{Conclusion}
\label{sec:conclusion}
In this paper we have described a particular statistical model and a Bayesian learning algorithm to fit it to data of high dimension, motivated by the source separation problem for Cosmic Microwave Background. The computational bottleneck of the implementation of this algorithm is to solve a linear system and we have shown 2 approximate solvers applied to this learning algorithm that make use of a Krylov subspace. We have shown that both approaches produce solutions at the size of Planck data (currently the most detailed all-sky data that we have) in a reasonable computation time with reasonable memory requirements. We show that the solution of this problem admits a representation as the solution of sparse-dense Sylvester equations.  For problems of the size generated by the \emph{real scientific problem from Planck satellite data}, our experiments indicate that standard approaches for sparse-dense problems become too costly to apply.  By instead using a modified version of the approach suggested in \cite[Section 5.2]{Palitta2018a}, our approach to solving the Sylvester equations has a significant advantage in computation time, taking only about one-sixth of the time of the conjugate gradient approach for larger data sizes.  This highlights the advantage of exploiting the Sylvester structure over treating the Kronecker-structured problem directly.
The almost-linear complexity of the algorithm implies that it will also be practical to use these approaches on any future data that is of larger size.

We note that a full implementation of the statistical solution also learns about the model parameters, such as those that define the $\vek A$ and $\vek Q$ matrices.  In this paper those have been assumed fixed, in order to focus on the issue of the linear system.  Since there are a small number of these parameters, the usual solution involves repeatedly solving the linear system for different combinations of  values of these parameters and building up a picture of which fit the data best. In Bayesian learning this can be done through the integrated nested Laplace approximation \cite{rue09}. Having to solve the linear system multiple times makes the need for an efficient solver even more important.  Thus, further improvements to accelerate the solution should be explored.  For example, subspace augmentation and recycling techniques (see, e.g., \cite{SoodhalterDeSturlerKilmer:2020:1} and references therein) would have the potential to accelerate the solution of multiple systems needed to learn the model parameters.

\section*{Acknowledgements}
The authors would like to thank the anonymous referees for their constructive critiques and suggestions, which greatly improved this manuscript.  The authors affirm they have no conflicts of interest.

%% file: appendix.tex
\appendix
\section{HEALPix}
\label{app:healpix}
HEALPix (Hierarchical Equal Area isoLatitude Pixelation) is a multi-resolution partition of the sphere into elements or pixels of equal area, although not the same shape.  It was originally used to provide a discretization of  all-sky observations done in connection with CMB but it is now used in many other fields. As of writing, an active website of resources is available at \url{https://healpix.sourceforge.io/}.

\begin{figure}
    \centering
	\includegraphics[scale=0.6]{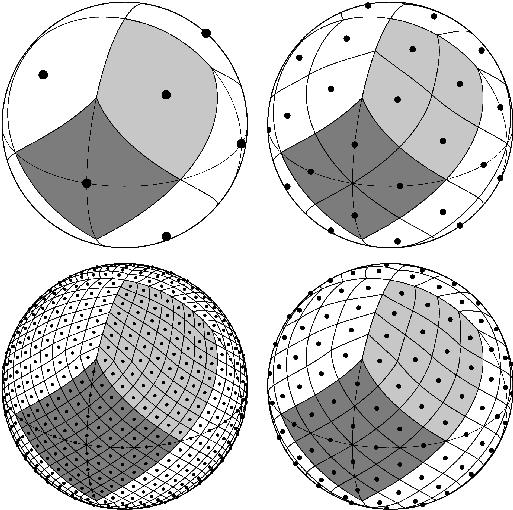}
    \caption{HEALPix partitions. Clockwise from top left: the base partition into 12 pixels, then level 1, 2 and 3 (48, 192 and 768 pixels). Taken from \url{https://healpix.sourceforge.io/}.}
    \label{fig:healpix}
\end{figure}

The base partition is into 12 elements. Successively more refined levels of the partition are made by dividing an element into 4 new elements (see Figure \ref{fig:healpix}). Taking the base partition to be level 0, the HEALPix partition at level $h$ therefore has $12 \times 4^h$ elements. For WMAP data, level 9 was used while for Planck a higher resolution of level 10 was possible, hence there are $12 \times 4^{10} = 12,582,912$ pixels in each image. 
HEALPix software supports data input and output, visualization and features like identifying neighbouring pixel indices and the indices of pixels in discs, triangles, polygons and strips.

Most analyses with data indexed using HEALPix will do separate analyses on each of the 12 base partition elements because functions such as identifying neighbouring pixels are only available within them.  That is what is done here. So our source separation on data at HEALPix level $h$ will consist of 12 independent analyses with $N = 4^h$.

\section{The matrix $\vek A$ in the CMB source separation problem}
\label{app:A}
The matrix $\vek A$ gives the coefficients of the assumed linear relationship between sources and observed brightnesses (Equation \ref{eq:model}). It is an $n \times m$ matrix and so typically small e.g.\ $n=9$ for Planck data and we might use $m=4$ or 5. The element $A_{kl}$ can be interpreted as the contribution of the $l$th source to the $k$th map. It is called the mixing matrix in the context of source separation.  CMB is indexed as the first source, so the first column of $A$ gives the contributions of CMB. In this work, 3 other sources are assumed, in order: synchrotron radiation, galactic dust and free-free emission.

The first thing to note is the well-known property of source separation solutions (and more generally factor analysis) that $\vek A$ and $\vek s$ are unique only up to an orthogonal transformation.  In the CMB application, uniqueness is gained by considering the separation in units of the CMB intensity, so that the first column of $\vek A$ is a vector of 1's.

For the other columns of $\vek A$, the physics of the sources are used to give expressions for the $A_{kl}$.  For synchrotron emission, a straight power law can be used:
% when considered in units of antenna temperature. Since CMB is best described in thermodynamic temperature, $a_l(\nu,\nu_0)$ becomes a scaled power law, multiplied by a conversion factor $c(\nu)$ from antenna to thermodynamic temperature 
\[
% a_{\mathrm{s}}( \nu, \nu_0 ) =
A_{k2} = c(\nu_k) \left(\frac{\nu_k}{\nu_0} \right)^{\kappa_s},
\]
for a reference frequency $\nu_0$ at which the source maps will be expressed (for Planck data this is usually the 4th frequency at 100 GHz) and $c(\nu) = (\mathrm{e}^{\psi} - 1 )^2/\psi^2 \mathrm{e}^{\psi}$, where $\psi = h \nu / k_B T_1$ with $h$ and $k_B$ are the Planck and Boltzmann constants and  $T_1 = 18.1\,\mathrm{K}$. The parameter $\kappa_s$ is known as the spectral index and is typically somewhere in the range (-2.3, -3.0).  Here we assume $\kappa_s = -2.65$.

For galactic dust, a modified black body spectrum is used, multiplied by a power law,
\[
% a_{\mathrm{d}}( \nu, \nu_0 ) 
A_{k3} = c(\nu_k) \frac{B(\nu_k)}{B(\nu_0)} \left(\frac{\nu_k}{\nu_0} \right)^{\kappa_d},
\] 
where $B(\nu_k) = \nu_k /\, [\, \exp( h \nu_k / k_B T_1) - 1\,]$. The spectral index for dust $\kappa_d$ is typically in the range $(1,2)$.  Here we assume $\kappa_d = 1.5$.

For free-free emission, a straight power law with known spectral parameter of $-2.14$ is used:
\[
% a_{\mathrm{ff}}(\nu,\nu_0) 
A_{k4} = c(\nu_k) \, \left(\frac{\nu_k}{\nu_0} \right)^{-2.14}.
\] 
For Planck data, where the 9 maps are for microwave frequencies at 30, 44, 70, 100, 143, 217, 353, 545 and 857 GHz, this yields
\[ \vek A = \left( \begin{array}{cccc} 
              1.000 & 24.314 &  0.181 & 13.158 \\
              1.000 &  8.817 &  0.315 &  5.801 \\
              1.000 &  2.581 &  0.612 &  2.151 \\
              1.000 &  1.006 &  1.006 &  1.006 \\
              1.000 &  0.392 &  1.630 &  0.471 \\
              1.000 &  0.132 &  2.783 &  0.196 \\
              1.000 &  0.038 &  4.931 &  0.072 \\
              1.000 &  0.013 &  7.704 &  0.032 \\
              1.000 &  0.005 & 11.337 &  0.015 \\
              \end{array} \right). \]       
We also need the vector $\vek T$ in the definition of $\vek C$.  For this application, the values are:
\begin{equation}\nonumber
\vek T = \mbox{diag}\{629881.6, 694444.4, 783146.7, 12755102.0, 30864197.5, 30864197.5, 30864197.5, 30864197.5, 30864197.5\}.
\end{equation}

\section{The intrinsic Gaussian Markov random field}
\label{app:igmrf}
A full derivation of the intrinsic Gaussian Markov random field (IGMRF) and its properties can be found in \cite[Chapter 3]{rue05}, of which below is a brief summary.

A set of random variables $X = (X_1,\ldots,X_N)$ is a \textit{Markov random field} if
\begin{itemize}
    \item Each variable $X_i$ has a neighbourhood set $\eta_i \subset \{1,\ldots,N\}$ that satisfies the properties: $i \notin \eta_i$ and $j \in \eta_i$ iff $i \in \eta_j$.  Intuitively, when the $X_i$ are values of pixels in a 2-dimensional image, $\eta_i$ is the indices of the neighbour pixels of pixel $i$. In this application using the HEALPix tesslation, the neighbours are the 4 pixels that share a common edge (as can be seen in Figure \ref{fig:healpix}).
    \item Each $X_i$ is independent of the other variables in the set given its neighbours. In other words $p(x_i \, | \, \{ x_j \, | \, j \neq i\}) = p(x_i \, | \, \{ x_j \, | \, j \in \eta_i \})$.
\end{itemize}
 Under quite general conditions, the joint probability distribution $p(x_1,\ldots,x_n)$ is uniquely defined by the set of conditional distributions $p(x_i \, | \, \{ x_j \, | \, j \in \eta_i \}), \; i=1,\ldots,N$, of each random variable given its neighbours. For the case of a 2-dimensional grid, the effect of these assumptions is that the values of nearby random variables on the grid are more correlated than those that are far away, and so realizations of the set of random variables are more spatially smooth than the case where each variable is independent.
 
 When these conditional distributions are Gaussian, the joint distribution is a multivariate Gaussian with some mean vector $\bmu$ and covariance matrix $\boldsymbol\Sigma$, and it is known as a \textit{Gaussian Markov random field} (GMRF).  The inverse of the covariance matrix, $\vek Q = \boldsymbol\Sigma^{-1}$, known as the precision matrix, has a nice structure since an off-diagonal element is zero if and only if the corresponding random variables are not neighbours:
 \[ j \notin 
 \{ i \} \cup \eta_i \: \Longleftrightarrow Q_{ij} = 0. \]
 In most cases, including ours, since the number of neighbours is small relative to $N$, $\vek Q$ will be sparse. % Often, permutation of the indices will yield $Q$ to be a band matrix with relatively small width, permitting efficient computation of its Cholesky decomposition that is sufficient to compute the multivariate Gaussian density and also for simulation.
 
 Finally, we get to the intrinsic GMRF. The key property of these models is that they are improper - they look like GMRFs but their joint probability distribution takes the form of a multivariate Gaussian but with a precision matrix that is not of full rank. Hence their joint distribution is not a properly defined probability density function.  This may seem to be a disqualifying property for a probability model but in fact (a) they model properties of the random variables that can be useful and (b) when used as a prior distribution, as in this application, the resulting posterior distribution is very often a properly defined probability distribution.
 
The vector of random variables $X = (X_1,\ldots,X_N)$ is an \textit{improper GMRF} of rank $N-k$ with parameters $(\bmu,\vek Q)$ if its density is of the form:
\[ p(x) = (2\pi)^{-0.5(N-k)} \, (|\widehat{\vek Q}|)^{0.5} \, \exp\left(-0.5(\vek x-\bmu)^T \vek Q (\vek x-\bmu) \right), \; x \in \mathbb{R}^N, \]
where $\vek Q$ is an $N \times N$ semi-positive definite matrix of rank $N-k > 0$ and $|\widehat{\vek Q}|$ is a generalised determinant, being the product of the $N-k$ non-zero eigenvalues of $\vek Q$.

The properties of such a construction are discussed in great detail in \cite{rue05}. A specific example of an IGMRF, that we use in this work, is where $\bmu = \vek 0$ and $\vek Q  = \varphi \vek D^T \vek D$, for $\varphi > 0$ and $\vek D$ an $N \times N$ matrix that encodes the neighbourhood structure,
\begin{equation}
D_{ij} \: = \: \begin{cases} 1, & \mbox{if } i \in \eta_j, \\ 0, & \mbox{otherwise}, \end{cases}
\label{eq:D}
\end{equation} 
giving $\vek Q$ a rank of $N-1$.  One can show that under this definition of $\vek Q$, the distribution of $X_i$ given its neighbours $\{ X_j \, | \, j \in \eta_i \}$ is  Gaussian with mean equal to the average of its neighbours. This makes a realisation of $X$ somewhat smooth, in that neighbouring values in the HEALPix lattice are more likely to take similar values than those that are further apart.  Another property of this construction is that there is no mean level specified --- although $\bmu = \vek 0$ this is not the 'mean' of $X$ as it is not a properly defined probability distribution.  The construct only gives properties of the relationship between neighbouring values in the lattice, not about their overall level.  This is useful if all that one wants to model for $X$ is this smooth spatial property, and not commit to a particular mean level a priori.  

That is the approach that we have taken in modelling the sources in our separation problem. Each source is modelled as an independent IGMRF, with the $l$th source having parameters $\vek{\mu}_l = 0$ and $\vek Q_l  = \varphi_l \vek D^T \vek D$, with $\vek D$ defined as in Equation \ref{eq:D}.  In this paper we have fixed $\varphi_l=1$ for each source. Stacking these independent sources into the vector $\vek S$ yields a Gaussian form with mean 0 and a block-diagonal precision matrix $\vek Q = \mbox{diag}(\vek Q_1,\ldots,\vek Q_m)$.